\documentclass[a4paper,12pt,draftcls,onecolumn]{IEEEtran}
\usepackage{amsmath}
\usepackage{amsfonts}
\usepackage{amssymb}
\usepackage{psfrag}
\usepackage{color}
\usepackage{eepic}
\usepackage{times}
\usepackage{cite}
\usepackage{multicol}
\usepackage{verbatim}
\usepackage{graphicx}
\usepackage{cases}
\usepackage{float}
\usepackage{pmat}
\usepackage{hyperref} 
\usepackage{eurosym}
\usepackage{acronym}
\usepackage[acronym,shortcuts]{glossaries}
\usepackage{cancel}
\usepackage{scalerel}
\usepackage{mathbbol}
\usepackage{tikz}
\usetikzlibrary{calc}
\usepackage[mathscr]{euscript}
\usepackage{wasysym}
\usepackage{subfigure}
\usepackage{color}
\usepackage{url}
\usepackage{arydshln}
\usepackage{algorithm}
\usepackage[noend]{algpseudocode}
\makeatletter
\def\BState{\State\hskip-\ALG@thistlm}
\makeatother

\DeclareSymbolFont{bbold}{U}{bbold}{m}{n}
\DeclareSymbolFontAlphabet{\mathbbold}{bbold}

\newcommand*\canc[1]{%
  \mathchoice
    {\scriptstyle#1}
    {\scriptstyle#1}
    {\scriptscriptstyle#1}
    {\scriptscriptstyle#1}
}

\newcommand*\Dcancelto[2][0]{%
  \kern9pt%
  \begin{tikzpicture}[baseline=(current bounding box.center).anchor=west]
    \node[anchor=east,inner sep=2pt] (a) {#2};
    \draw[->] ($(a.north west)+(1pt,-2pt)$) -- ($(a.south east)+(0pt,2pt)$) node at ($(a.south east)+(4pt,1pt)$) {$\canc{#1}$};
\end{tikzpicture}
}

\renewcommand{\d}{{\rm d}}
\newcommand{\tr}{{\rm tr}}
\newcommand{\sol}{{\rm sol}}

\newacronym{CRLB}{CRLB}{Cram{\'e}r-Rao lower bound}
\newacronym{NLOS}{NLOS}{non line-of-sight}
\newacronym{LOS}{LOS}{line-of-sight}
\newacronym{WSN}{WSN}{wireless sensor network}
\newacronym{ToA}{ToA}{time of arrival}
\newacronym{AoA}{AoA}{angle of arrival}
\newacronym{DoA}{DoA}{direction of arrival}
\newacronym{PDoA}{PDoA}{phase-difference of arrival}
\newacronym{SNR}{SNR}{Signal-to-Noise Ratio}
\newacronym{TWR}{TWR}{two-way ranging}
\newacronym{CWRR}{CWRR}{continuous wave radar ranging}
\newacronym{RSSI}{RSSI}{receive signal strength indicator}
\newacronym{OFDM}{OFDM}{orthogonal frequency division multiplexing}
\newacronym{MUSIC}{MUSIC}{Multiple Signal Classification}
 
\begin{document}
\title{Superresolution Multipoint Ranging with\\[-1ex] Optimized Sampling via Orthogonally\\[-1ex] Designed Golomb Rulers}

\author{Omotayo Oshiga,~\IEEEmembership{Student Member,~IEEE,} Stefano Severi,~\IEEEmembership{Member,~IEEE,}\\ and 
        Giuseppe~Abreu,~\IEEEmembership{Senior Member,~IEEE}
\thanks{O.~Oshiga, S.~Severi and G.~Abreu are with the School of Engineering and Sciences, Jacobs University Bremen,
Campus Ring 1, 28759 Bremen, Germany. E-mails:{\tt [o.oshiga,s.severi,g.abreu]@jacobs-university.de}}
\thanks{Parts of this material will be presented in the Wireless Communication and Networking Conference 2014 \cite{Oshiga2014}.}\vspace{-3ex}
}

\markboth{IEEE Transactions on Wireless Communications,~Vol.~X, No.~X, XXXX~201X}%
{Abreu: Bla bla bla}

\maketitle

\vspace{-6ex}
\begin{abstract}
\vspace{-2ex}
We consider the problem of performing ranging measurements between a source and multiple receivers efficiently and accurately, as required by distance-based wireless localization systems.
To this end, a new multipoint ranging algorithm is proposed, which is obtained by adapting superresolution techniques to the ranging problem, using for the sake of illustration the specific cases of \ac{ToA} and \ac{PDoA}, unified under the same mathematical framework.
The algorithm handles multipoint ranging in an efficient manner by employing an orthogonalized non-uniform sampling scheme optimised via Golomb rulers.
Since the approach requires the design of mutually orthogonal sets of Golomb rulers with equivalent properties -- a problem that founds no solution in current literature -- a new genetic algorithm to accomplish this task is presented, which is also found to outperform the best known alternative when used to generate a single ruler.
Finally, a \acp{CRLB} analysis of the overall optimised multipoint ranging solution is performed, which together with a comparison against simulation results validates the proposed techniques.
\end{abstract}

\vspace{-3ex}

\section{Introduction}
\label{Sec:Introduction}\glsresetall
\vspace{-1ex}

%
%
Wireless localization is a fairly mature area of research, with a vast literature \cite{Macagnano2012, Yanying2009, Hui2007}.
It is therefore paradoxical that despite the formidable effort put into the problem, wireless positioning is still shy of its potential as a truly ubiquitous technology \cite{Harrop2012, Wirola2010, Hui2007, Harrop2013}.

Ubiquity requires the technology to be available in every environment, and it is well-known that wireless localization systems are still inaccurate and unreliable in places such as urban canopies and indoors, which are characterized by rich multipath and scarcity of \ac{LOS} conditions.
Furthermore, compared to the quality and omnipresence of satellite- and cellular-based systems in open outdoor spaces, indoor positioning solutions \cite{Convergence2012, Dresden2013, Medina2013} are still relatively fragile, under-deployed and unconsolidated.


One explanation for this discrepancy is that literature has provided a large number of building blocks to solve parts of the problem, but which for a reason or another still do not come together harmoniously to provide a comprehensive solutions.

To qualify the latter statement, consider the specific case of \ac{AoA} or \ac{DoA} positioning.
A good number of AoA-based localization algorithms \cite{Niculescu2003, Rong2006, Kegen2007, Azzouzi2011}, and an even wider body of literature on AoA estimation \cite{Schmidt1986, Barabell1983, FuLi1993, Gershman1999, tuncer2009} exists.
Of particular relevance is the fact that simultaneous estimation of the AoA of multiple signals/sources is relatively easy to perform, which is of fundamental importance to reduce latency in indoor applications where the concentration of users is typically large.

Yet, AoA-based indoor positioning is not common today because: $a$) AoA-based localization algorithms are highly susceptible to \ac{NLOS} conditions, such that accurate and robust AoA input is needed; and $b$) accurate and robust AoA estimation requires expensive multi-antenna systems and high computational capabilities, which are incompatible with typical indoor requirements of small, low-cost, low-power devices \cite{Wirola2010, Harrop2013}.

On the other extreme of the technological spectrum are proximity-based (in particular RFID) approaches \cite{Hui2007, Ristic2006, Shigeng2010}, which do satisfy the latter requirements, but at the expense of accuracy, and therefore also failed to penetrate the general market.

The limitations of the AoA- and proximity-based approaches partially explain the predominance of range-based indoor localization systems proposed both by academia \cite{Zhang2005, Hui2007, Shouhong2009, Yanying2009, Moragrega2010, Macagnano2012, JunlinYan2013}.
Indeed, various accurate and robust distance-based localization algorithms exist, and distance estimates are relatively inexpensive to obtain from radio signals -- via \ac{RSSI}, \ac{ToA} or \ac{PDoA} methods -- without requiring multiple antennas or significant additional RF circuitry.
But again the deployment of this technology is short of its potential, which arguably is a result of the fact that since ranging quality is severely degraded by interference, range-based positioning systems are required to carefully schedule the collection of ranging information, leading to low refreshing ratios and high communication costs.

The above rationale points to a curious predicament.
On the one hand, many excellent multipoint \ac{AoA} estimation algorithms exist \cite{Schmidt1986, Barabell1983, FuLi1993, Gershman1999, tuncer2009}, which however are not typically utilised for indoor positioning as multi-antenna systems are too expensive.
On the other hand, many excellent distance-based localization algorithms exist \cite{Zhang2005, Hui2007, Shouhong2009, Yanying2009, Moragrega2010, Macagnano2012, JunlinYan2013}, which however can only be effectively employed for indoor positioning if ranging information can be collected efficiently from multiple sources so as to reduce latency.


The work presented in this article is a proposal to solve the aforementioned impasse.
Specifically, we offer a solution to the multipoint ranging problem based on the same superresolution techniques typically used for AoA estimation.
As shall be explained, however, in this context the ability to handle multiple sources when employing superresolution methods does not stem from the separability of signals through the eigen-properties of mixed covariance matrices, but rather by a robustness to sampling sparsity which interestingly is not always enjoyed by such methods in the multi-antenna setting.
The feature suggests that the collection of input data can be optimized by designing such sampling sparsity according to Golomb rulers \cite{Dewdney1985, Rankin1993, Soliday1995, CottaCONSTRAINTS2007}, which however must maintain mutual orthogonality.
The latter is achieved by a new genetic algorithm -- designed under the inspiration of the behaviour of prides of lions -- which enables the construction of multiple orthogonal and equivalent\footnote{Equivalence will be defined more rigorously according to two different criteria.} Golumb rulers.

The performance of the new algorithm to construct Golomb rulers is compared against the state of the art, and shown thereby to outperform all alternatives we could find.
Furthermore, an original \acp{CRLB} analysis of the new strategy is performed, which indicates that in addition to the advantage of enabling simultaneous multipoint ranging, the overall solution achieves remarkable gain in accuracy over current methods.

In summary, our contributions are as follows:
\begin{itemize}
\item[1)] A new multipoint ranging algorithm obtained by adapting superresolution techniques for \ac{ToA} \cite{RainerHach2005, Baba2011, Myungkyun2010, JianXing2007} and \ac{PDoA} \cite{Scherhaufl2013, Povalac2011, Ahmad2006} ranging, under a unified mathematical framework;
\item[2)] A new genetic algorithm that outperforms the best known alternative and enables the construction of multiple orthogonal sets of Golomb rulers of equivalent properties;
\item[3)] A complete \acp{CRLB} analysis of the resulting method, which validates its advantages.
\end{itemize}

%
%
%
%

\section{Super-resolution {ToA} and {PDoA} Ranging}
\label{Sec:Prelim}

%
%

There are three basic methods to estimate the distance between a pair of wireless devices using their signals: \ac{RSSI}, \ac{ToA} and \ac{PDoA}.
Amongst these alternatives, \ac{RSSI}-ranging is known to be the least accurate and least robust \cite{Chandrasekaran2009, Elnahrawy2004}.
In fact, after some early attention due mostly to its inherent low-power potential \cite{Qianqian2011,Keping2013}, \ac{RSSI}-ranging has since lost appeal thanks to the emergence of low-power physical layer standards such as 802.15.4g \cite{IEEE802.11-2011} and 802.11ac \cite{Szulakiewicz2012}, which facilitate the implementation of low-power \ac{ToA} and \ac{PDoA} ranging mechanisms.
In light of the above, we shall focus hereafter on the latter two forms of ranging.

\subsection{ToA-based Two-Way Ranging}
\label{Subsubsec:Twr}

Consider the problem of estimating the distance $d$ between a reference node (anchor) $A$ and a target node $T$ based on ToA measurements.
Using the standard two-way ranging technique \cite{RainerHach2005, Baba2011, Myungkyun2010, JianXing2007}, and assuming that the procedure is executed not a single but multiple times, the $k$-th distance estimate $\hat{d}_k$ of $d$ is computed by
\begin{equation}
\label{Eq:distancetoa}
\hat{d}_k = \Big[\big(\tau_{_{\scriptstyle\textup{RX}:k}} - \tau_{_{\scriptstyle\textup{TX}:k}}\big) - k\cdot\tau_{_{\scriptstyle T}}\Big]\cdot\dfrac{c}{2}
\end{equation}
where $c$ is the speed of light; $\tau_{_{\scriptstyle\textup{TX}:k}}$ and $\tau_{_{\scriptstyle\textup{RX}:k}}$ are respectively the time stamps of the $k$-th packet at transmission and reception back at the anchor; and  $\tau_{_{\scriptstyle T}}$ is a fixed and known waiting period observed by the target, for reasons that are beyond\footnote{For instance, $\tau_{_{\scriptstyle T}}$ may be imposed by the frame structure of the underlying communication system.} the ranging process itself.

Since $\tau_{_{\scriptstyle T}}$ is known \emph{a priori} by the anchor, it serves no mathematical purpose and therefore can be assumed to be zero\footnote{Strictly speaking, $\tau_{_{\scriptstyle T}}$ could also be considered a source of ranging errors, since it is subject to jitter (imperfect time-keeping). In practice, however, jitter errors are several orders of magnitude below the timing errors involved in measuring $\tau_{_{\scriptstyle\textup{RX}:k}}$, and therefore can be effectively ignored.} without loss of generality (w.l.g.).
Similarly, before the $k$-th ranging cycle the anchor may in practice hold for a (possibly unequal) waiting period $\tau_{_{\scriptstyle\! A:(k-1)}}$, which however can also be normalized to zero, wlg.

\begin{figure}[H]
\centering
\includegraphics[width=0.6\columnwidth]{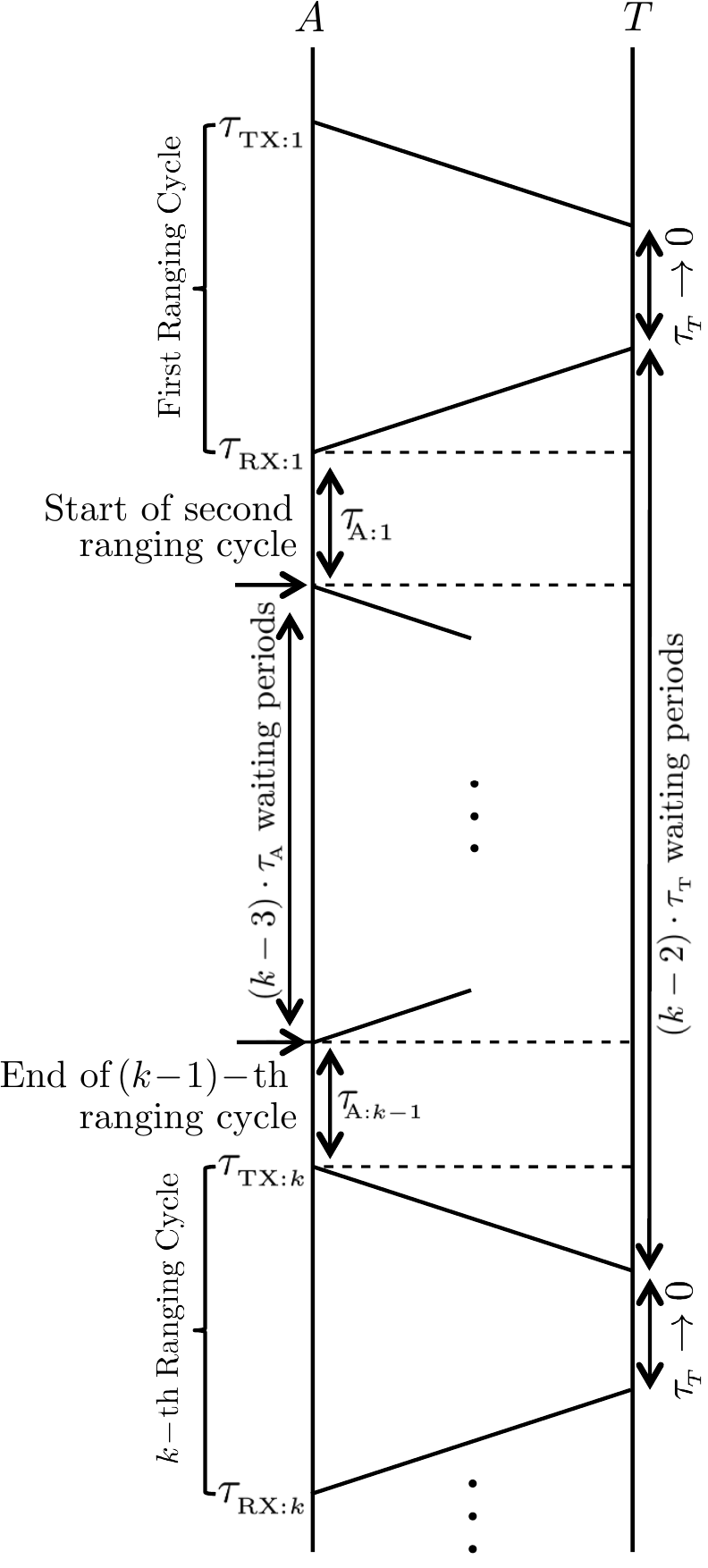}
\vspace{-3ex}
\caption{Illustration of non-uniform TWR scheme. Multipoint-point ranging can be performed by intercalating different sources in different orthogonal (non-overlapping) slots (cycles).}
\label{Fig:ToA_Top}
\end{figure}

Referring to Figure \ref{Fig:ToA_Top}, and considering the latter assumptions on $\tau_{_{\scriptstyle T}}$ and $\tau_{_{\scriptstyle\! A:i}}$ for $i=\{1,\cdots,k-1\}$, equation \eqref{Eq:distancetoa} can then be rewritten as
\begin{equation}
\label{Eq:Mul_two_way_rang}
\hat{d}_k = \left[\underbrace{\big(\tau_{_{\scriptstyle\textup{RX}:k}} -\tau_{_{\scriptstyle\textup{TX}:1}}\big)}_{\Delta \tau_k } -  k\cdot\!\!\!\!\!\!\Dcancelto[0]{$\tau_{_{\scriptstyle T}}$} - \sum_{i=1}^{k-1}\!\!\!\!\Dcancelto[0]{$\tau_{_{\scriptstyle\! A:i}}$}\!\right]\cdot\dfrac{c}{2 k} \equiv \Delta \tau_k \cdot\dfrac{c}{2 k}.
\end{equation}

One way to interpret the model described by equation \eqref{Eq:Mul_two_way_rang} is that in a \ac{ToA}-based \ac{TWR} scheme with multiple ranging cycles, the time-difference measurement $\Delta \tau_k$ obtained at the $k$-th cycle has a linear functional relationship with the cycle index $k$, with the proportionality factor determined by the distance $d$ between the target and the anchor, $i.e.$,
\begin{equation}
\label{Eq:TimeDifferenceDistance}
\Delta \tau_k  = \omega_d k, \quad \mbox{with} \quad \omega_d = \frac{2d}{c}.
\end{equation}

The convenience of this interpretation of \ac{ToA}-based \ac{TWR} will soon become evident.

\subsection{{PDoA}-based Continuous Wave Radar Ranging}
\label{Subsubsec:DFRR}
\begin{figure}[H]
\centering
\includegraphics[width=0.6\columnwidth]{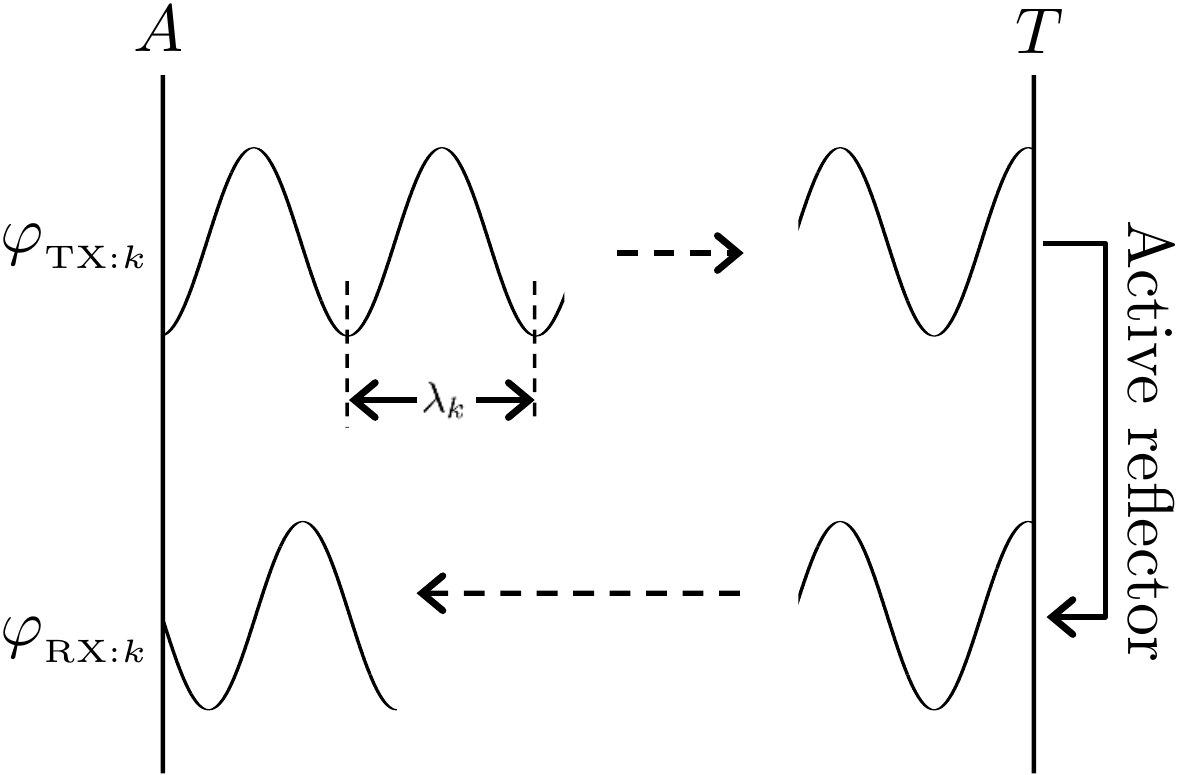}
\vspace{-3ex}
\caption{Illustration of PDoA ranging mechanism for a single frequency. Multipoint-point ranging can be performed by allocating different sources to different orthogonal carriers.}
\label{Fig:PDoA_Top} 
\vspace{7ex}
\includegraphics[height=0.4\textheight]{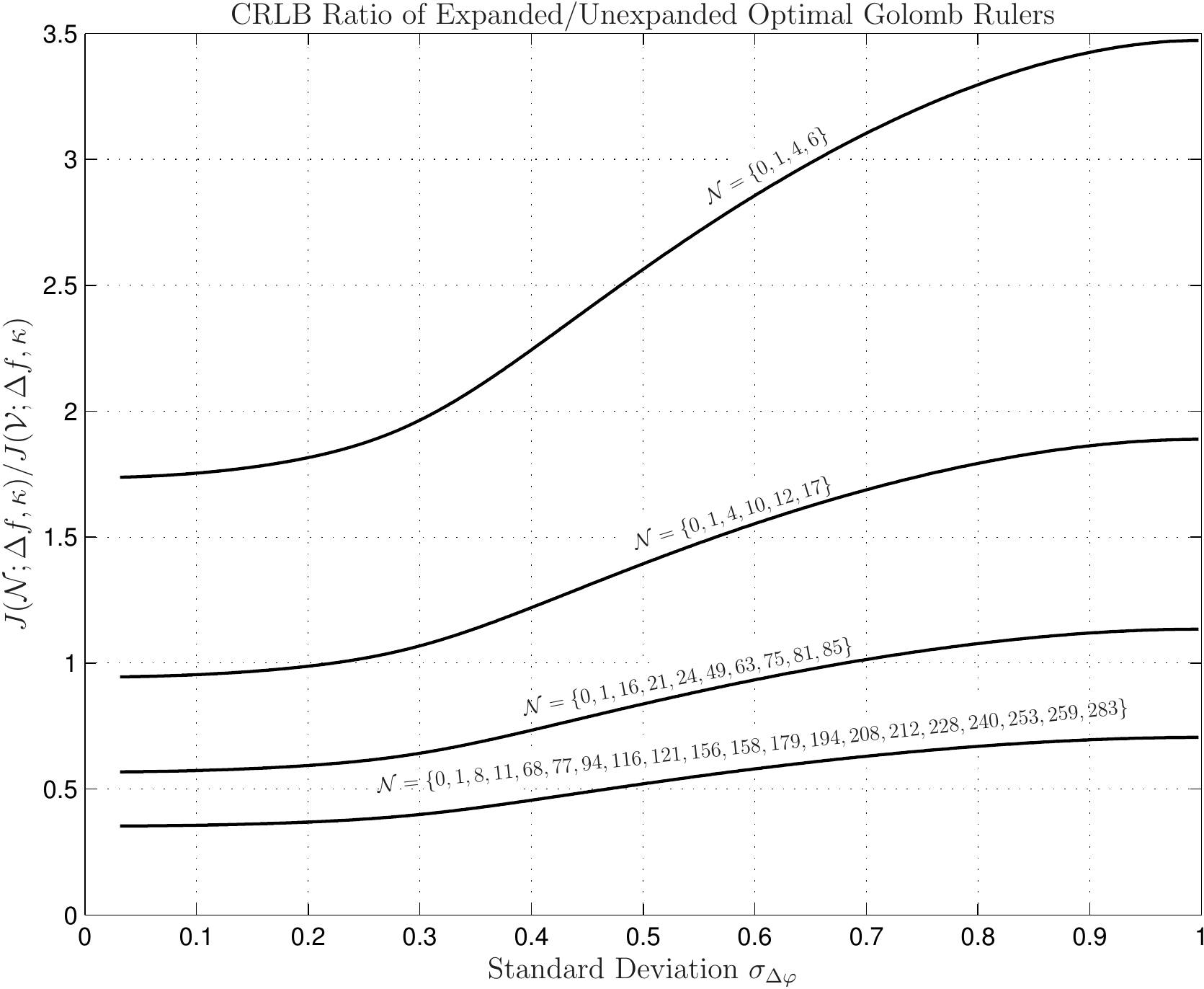}
\vspace{-3ex}
\caption{Evolution of CRLB ratio $J(\mathcal{N};\Delta f,\kappa)/J(\mathcal{V};\Delta f,\kappa)$ as a function of the phase error standard deviation $\sigma_{\Delta\varphi}$, associated with different rulers $\mathcal{N}$.}
\label{Fig:crlb_kappa} 
\end{figure}

Consider the problem of estimating the distance $d$ between a reference node (anchor) $A$ and a target node $T$ based on the phases of the signals exchanged between the devices.
One possible mechanism, as illustrated in figure \ref{Fig:PDoA_Top}, is that the anchor $A$ emits a continuous sinusoidal wave of frequency $f$ with a known phase $\varphi_{_{\textup{TX}}}$ and the
target $T$ acts as an active reflector, such that $A$ can measure the phase $\varphi_{_{\textup{RX}}}$ of the returned signal \cite{Scherhaufl2013, Povalac2011, Ahmad2006}.
In this case, the roundtrip distance $2d$ and the phases $\varphi_{_{\textup{TX}}}$ and $\varphi_{_{\textup{RX}}}$ are related by 
\begin{equation}
\label{Eq:CWRR}
\varphi = \varphi_{_{\textup{RX}}} - \varphi_{_{\textup{TX}}} = \dfrac{4 \pi d}{c} f - 2\pi N,
\end{equation}
where $N$ is the integer number of complete cycles of the sinusoidal over the distance $2d$.

Obviously the distance $d$ cannot be estimated directly based on equation \eqref{Eq:CWRR} since the quantity $N$ is unknown.
However, taking the derivative of equation \eqref{Eq:CWRR} with respect to $f$ one obtains
\begin{equation}
\label{Eq:CWRRDerivative}
\dfrac{\d \varphi}{\d f} = \dfrac{4 \pi d}{c}.
\end{equation}

Let there be a set of equi-spaced frequencies $\mathbb{F}=\{f_0,\cdots,f_K\}$ such that $\Delta f = f_{k+1} - f_k$ for all $0\leq k < K$,
and assume the roundtrip phases $\varphi_k$ for all $f_k$ are measured.
Then, thanks to the linear relationship between $f$ and $d$ described by equation \eqref{Eq:CWRRDerivative}, it follows that 
\begin{equation}
\label{Eq:PhaseDifferenceDistance}
\Delta \varphi_k  = \omega_d k, \quad \mbox{with} \quad \omega_d = \frac{4\pi\Delta f d }{c},
\end{equation}
where $\Delta \varphi_k \triangleq \varphi_k - \varphi_0$ for all $1\leq k < K$.

Comparing equations \eqref{Eq:TimeDifferenceDistance} and \eqref{Eq:PhaseDifferenceDistance}, we conclude that both the \ac{ToA}-based \ac{TWR} and the \ac{PDoA}-based \ac{CWRR} methods are mathematically equivalent, in the sense that the measured quantities, respectively $\Delta \tau_k$ and $\Delta \varphi_k$, have a linear relationship with a counter $k$, governed by a slope coefficient $\omega_d$ that is directly and unequivocally related to the desired information $d$.

In light of the models described above, we shall consider for simplicity that we are able to measure  quantities $\Delta_k$, such that
\begin{equation}
\label{Eq:GenericRanging}
\Delta_k  = \omega_d \cdot k,
\end{equation}
where $\omega_d$ is a coefficient with a constant relationship with $d$.

Notice that trivially due to the linearity of this relationship, we have, for any pair of integers $(k,q)$,  with $k > q$,
\begin{equation}
\label{Eq:GenericRangingDifference}
\Delta_k - \Delta_q  = \omega_d \cdot (k - q) = \Delta_{k - q}.
\end{equation}

This simple property has a remarkable consequence.
Indeed, consider an ascending sequence of non-negative integers $\mathcal{N}= \{n_1,\cdots,n_K\}$ and the associated set of input measurements $\scaleobj{1.3}{\mathbbold{\Delta}}_\mathcal{N} = \{\Delta_{n_1},\cdots,\Delta_{n_K}\}$.
By virtue of equation \eqref{Eq:GenericRangingDifference}, the set $\scaleobj{1.3}{\mathbbold{\Delta}}_\mathcal{N}$ can be expanded into $\scaleobj{1.3}{\mathbbold{\Delta}}_\mathcal{V} =
\{\Delta_{n_2}\!-\Delta_{n_1},\cdots,\Delta_{n_K}\!-\Delta_{n_1},\cdots,\Delta_{n_K}\!-\Delta_{n_{K-1}}\} = \{\Delta_{n_2-n_1},\cdots,\Delta_{n_K - n_{K-1}}\} = \{\Delta_{\nu_1},\cdots,\Delta_{\nu_M}\}$, where the cardinality $M$ of $\scaleobj{1.3}{\mathbbold{\Delta}}_\mathcal{V}$ is obviously upper bounded by $M \leq K\frac{K-1}{2}$.

Other then the much larger cardinality, the sequences $\scaleobj{1.3}{\mathbbold{\Delta}}_\mathcal{V}$ and $\scaleobj{1.3}{\mathbbold{\Delta}}_\mathcal{N}$ have, as far as the purpose of distance estimation is concerned, fundamentally the same nature since both carry samples of the quantities $\Delta_k$.
In other words, the model described in subsection \ref{Subsubsec:DFRR} allows for large input sets of cardinality $N$ to be obtained from a significantly smaller number $K$ of actual measurements, by carefully designing the feedback intervals or the carrier frequencies required to perform ranging estimates.
Furthermore, the linearity between the measured quantities $\Delta_k$ and the corresponding indexes $k$ is so that such design can be considered directly in terms of the relationship between the integer sequences $\mathcal{N} \to \mathcal{V}$.

Sparse sequences $\mathcal{N}$ that generate optimally expanded equivalents $\mathcal{V}$ are known as \emph{Golomb rulers} and their design under the constraints of our problem is the subject in Section \ref{Sec:Gol_Sr_R}.
Here, however, let us proceed by demonstrating how the aforementioned model enables the straightforward application of superresolution algorithms for \ac{ToA} and  \ac{PDoA} ranging.

%
%
%
%


\vspace{-2ex}
\subsection{Multi-point Ranging via Super-resolution Algorithms}
\label{Subsec:S_PDB_R}

Straightforwardly, assume that a set of input measurements $\scaleobj{1.3}{\mathbbold{\Delta}}_\mathcal{N}$ is collected, from which the associated expanded set $\scaleobj{1.3}{\mathbbold{\Delta}}_\mathcal{V}$ is constructed and consider the corresponding complex vector
\vspace{-3ex}
\begin{equation}
\label{Eq:Steerin_vect}
\mathbf{x} = [e^{j\Delta_{\nu_1}},e^{j\Delta_{\nu_2}},\cdots,e^{j\Delta_{\nu_M}}]^\textup{T} \equiv [e^{j\omega_d},e^{j{\nu_2}\omega_d},\cdots,e^{j{\nu_M}\omega_d}]^\textup{T},
\vspace{-1ex}
\end{equation}
where $^\textup{T}$ denotes transposition and we have normalised ${\nu_1}=1$, without loss of generality.

One can immediately recognize from equation \eqref{Eq:Steerin_vect} the similarity between the vector $\mathbf{x}$ and the steering vector of a linear antenna array \cite{Schmidt1986, Barabell1983, tuncer2009}, with inter-element spacings governed by $\scaleobj{1.3}{\mathbbold{\Delta}}_\mathcal{V}$.
An estimate of the parameter of interest $\omega_d$ can therefore be recovered from the covariance matrix $\mathbf{R}_\mathbf{x} \triangleq \mathbb{E}[\mathbf{x}\cdot\mathbf{x}^\textup{H}]$.
Specifically, under the assumption that each measurement $\Delta_{\nu_m}$ is subject to independent and identically distributed (iid) white noise with variance $\sigma^2$, the covariance matrix $\mathbf{R}_\mathbf{x}$ can be eigen-decomposed to
\vspace{-1ex}
\begin{equation}
\mathbf{R}_\mathbf{x} = \mathbf{U}\cdot\boldsymbol{\Lambda}\cdot\mathbf{U}^\textup{H},
\vspace{-1ex}
\end{equation}
with
\vspace{-1ex}
\begin{equation}
\mathbf{U} = 
\begin{pmat}[{|c}]
\mathbf{u}_\mathbf{x} & \mathbf{U}_0,\cr
\end{pmat}
\quad \mbox{and} \quad
\boldsymbol{\Lambda} =
\begin{pmat}[{|c}]
1 + \sigma^2 & \boldsymbol{0}\cr
\-
\boldsymbol{0} & \sigma^2 \mathbf{I}\cr
\end{pmat},
\end{equation}
where $\mathbf{U}_0$ is the $K$-by-$(K\!-\!1)$ null-space of $\mathbf{R}_\mathbf{x}$.

Given the above properties, many superresolution algorithms can be employed to obtain ranging estimates from \ac{ToA} and \ac{PDoA} measurements \cite{Rubsamen2009, Xiong2012, Abdalla2013,Faye2013}.
Since our focus in this article is to demonstrate such possibility, discuss the resulting opportunities to optimize resources, and analyze the corresponding implications on the achievable ranging accuracies, we shall limit ourselves to two explicit classical examples, for the sake of clarity.

One way to obtain an estimate $\hat{\omega}_d$ of $\omega_d$ is via the classic spectral \ac{MUSIC} algorithm \cite{Schmidt1986, Gershman1999, Stoica1990}, where a search for the smallest vector projection onto the noise subspace of $\mathbf{R}_\mathbf{x}$ is conducted, namely
\begin{equation}
\label{Eq:MUSIC}
\hat{\omega}_d = {\rm arg}\max\limits_{\omega_d} \frac{1}{\|\mathbf{e}^\textup{H}\cdot \mathbf{U}_0\|^2}
\quad \mbox{with} \quad \mathbf{e}\triangleq [e^{j\omega_d},e^{j{\nu_2}\omega_d},\cdots,e^{j{\nu_M}\omega_d}]^\textup{T}.
\end{equation}

Alternatively, $\hat{\omega}_d$ can be obtained using the root \ac{MUSIC} algorithm \cite{Barabell1983, tuncer2009, ElKassis2010}, which makes use of the fact that the projection square norm $\|\mathbf{e}^\textup{H}\cdot \mathbf{U}_0\|^2$ defines an equivalent polynomial in $\mathbb{C}$ with coefficients fully determined by the Grammian matrix of the null subspace of $\mathbf{R}_\mathbf{x}$.
Specifically, define the auxiliary variable $z \triangleq e^{j\omega}$ such that $\mathbf{e} = [z,z^{\nu_2},\cdots,z^{\nu_M}]^\textup{T}$, and the two zero-padded vectors 
$\mathbf{e}_{_\textup{L}} = [z^{-1},0,\cdots,0,z^{-\nu_2},0,\cdots,0,z^{-\nu_3},0,\cdots,\cdots,0,z^{-\nu_M}]$ and 
$\mathbf{e}_{_\textup{R}} = [z,0,\cdots,0,z^{\nu_2},0,\cdots,0,z^{\nu_3},0,\cdots,\cdots,0,z^{\nu_M}]$.
Then we may write
\begin{eqnarray}
\label{Eq:MUSICPolynomial}
P(z) \hspace{-3ex}&& = \|\mathbf{e}^\textup{H}\cdot \mathbf{U}_0\|^2 = \mathbf{e}_{_\textup{L}}\cdot \mathbf{G}\cdot \mathbf{e}_{_\textup{R}}^\textup{T} \equiv \sum_{\nu=0}^{2\nu_M-2} \tr(\mathbf{G};\nu)\cdot z^{\nu},
\end{eqnarray}
where the last equivalence sign alludes to the multiplication by $z^{\nu_M}$ required to take the algebraic function into a polynomial; $\mathbf{G}$ is a Gramian matrix constructed by zero-padding the matrix $\mathbf{U}_0\cdot\mathbf{U}_0^\textup{H}$, such that the $(m,\ell)$-th element of $\mathbf{U}_0\cdot\mathbf{U}_0^\textup{H}$ is the $(\nu_m,\nu_\ell)$-th element of $\mathbf{G}$; and $\tr(\mathbf{G};\nu)$ denotes the $m$-th trace of the matrix $\mathbf{G}$ -- $i.e.$, the sum of the $k$-th diagonal of $\mathbf{M}$, counting from the the bottom-left to the upper-right corner.

The estimate $\hat{\omega}_d$ can then be obtained by finding the only unit-norm root of $P(z)$, $i.e.$,
\begin{equation}
\label{Eq:RootMUSIC}
\hat{\omega}_d =  {\rm arg}\, \sol\, \Big\{P(z) = 0 \;\Big|\; |z| = 1\Big\}.
\end{equation}

Whatever the specific method used to extract the distance information (embedded in $\hat{\omega}_d$) from the vectors constructed as shown in equation \eqref{Eq:Steerin_vect}, the following properties apply to the superresolution algorithms described above.

\begin{itemize}
\item Superposibility: Thanks to the expansions $\mathcal{N} \to \mathcal{V}$, measurement intervals/frequencies corresponding to multiple sources can be superposed without harm.
To exemplify, consider the case of two sources $A$ and $B$ and the measurements from both sources be collected continuously according to the sequence $\mathcal{N} = \{1,3,4,5,6,7,8,10\}$, but such that the sources $A$ and $B$ are only active according to the orthogonal sequences $\mathcal{N}_A = \{1,3,6,7\}$ and
$\mathcal{N}_B = \{4,5,8,10\}$.
The samples in $\mathcal{N}_A$ can, however, be transformed into the sequence $\mathcal{V}_A = \{3-1,6-1,7-1,6-3,7-3,7-6\} \equiv \{1,2,3,4,5,6\}$, which contains $6$ samples.
Furthermore and likewise, $\mathcal{N}_B \to \mathcal{V}_B = \{5-4,8-4,10-4,8-5,10-5,10-8\}\equiv \{1,2,3,4,5,6\}$.
In other words, out of only 8 jointly collected samples, 6 \ac{ToA} or \ac{PDoA} (equivalent) measurements from each source are obtained, without interference.

\item Unambiguity: In the case of \ac{AoA} estimation using antenna arrays, the elements of the steering vectors are complex numbers whose arguments are \emph{periodic functions} of the desired parameter, which in turn gives rise to aliasing (ambiguity) of multiple parameter values that lead to the same set of measurements \cite{Song2002, Byungwoo2004, Keller2006, Tayem2012}.
In contrast, in the context hereby the quantities $\Delta_k$ are \emph{linear functions} of the desired parameter $d$, such that no such ambiguity occurs.
\item Separability: Thanks to both properties above, superresolution ranging can be carried without interference using orthogonal non-uniform sample vectors, each processed by a separate estimator.
Consequently, issues such as correlation amongst multiple signals, which commonly affect superresolution algorithms \cite{Chongying2007, Ariananda2012, Weiziu1991, Wang2002, Tayem2012}, do not exist in the context hereby.
In other words, the application of superresolution algorithms to multipoint ranging are more closely related to Pisarenko's original harmonic decomposition algorithm \cite{Pisarenko1973}, than to derivative methods such as \ac{MUSIC}.
\end{itemize} 


\section{Optimization of ToA and PDoA Range Sampling via Golomb Rulers}
\label{Sec:Gol_Sr_R}

Under the mathematical model described in Section \ref{Sec:Prelim}, the optimization of ranging resources amounts to allocating ranging cycles or frequency pairs to multiple sources, respectively, which is  directly related to that of designing Golomb rulers \cite{Rankin1993}.

Golomb rulers are sets of integer numbers that generate, by means of the difference amongst their elements, larger sets of integers, without repetition.
The problem was first studied independently by Sidon \cite{Sidon1932} and Babcock \cite{Babcock1953}, but these special sets are named after Solomon W. Golomb \cite{Golomb1977} as he was the first to popularize their application in engineering.
Before we discuss the design of Golomb rulers for the specific application of interest, it will prove useful to briefly review some of their basic characteristics and features.

\vspace{-3ex}
\subsection{Basic Characteristics and Features of Golomb Rulers}
\label{Subsec:Gol_Bas}

Consider a set of ordered, non-negative integer numbers $\mathcal{N} = \{n_1, n_2, \cdots,n_K\}$, with $n_1 = 0$ and $n_K = N$, wlg\footnote{Since Golomb rulers are invariant to translation, we consider without loss of generality, that the first element is $0$ and the last is $N$. That is slightly different from the representation adopted in subsection \ref{Subsec:S_PDB_R}, but will prove convenient hereafter.}.
This set has cardinality (or \emph{order}) $K$, and it will prove convenient to define the \emph{length} of
the set by its largest element $N$.

Next, consider the corresponding set $\mathcal{V}$ of all possible pairwise differences
\begin{equation}
\label{Eq:Golomb}
\nu_{k\ell} = n_k - n_\ell \quad (1 \leq \ell < k \leq K).
\end{equation}

If the differences $\nu_{k\ell}$ are such that $\nu_{k\ell} = \nu_{pq}$ if and only if (iff) $k=p$ and $\ell=q$, then the set
$\mathcal{N}$ is known as a \emph{Golomb ruler}.
Such sets are thought of as \emph{rulers}, as their elements can be understood as \emph{marks} of a ruler, which can thus \emph{measure} only the lengths indicated by any pair of marks.
In analogy to the latter, we henceforth refer to the set $\mathcal{V}$ as the \emph{measures} set.

It follows from the definition that the number of distinct lengths that can be measured by a Golomb ruler -- in other words, the order of $\mathcal{V}$ -- is equal to $K\frac{K-1}{2}$.
The first key feature of a Golomb ruler is therefore that if $\mathcal{N}$ has order $K$, then $\mathcal{V}$ has order $K\frac{K-1}{2}$.

A simple example of a Golomb ruler is $\mathcal{N} = \{0,1,4,6\}$, which generates the Measures $\mathcal{V}=\{1,2,3,4,5,6\}$.
In this particular example, $\mathcal{V}$ is \emph{complete}, as it contains all positive integers up to its length, so that the Golomb of order $4$ is said to be \emph{perfect}.
In other words, a perfect ruler allows for \emph{all lengths} to be measured, up to the length of the ruler itself.

Unfortunately, no perfect Golomb ruler exists \cite{Dewdney1985} for $K > 4$.
It is therefore typical to focus on designing rulers that retain another feature of the order-4 Golomb ruler, namely, its compactness or \emph{optimality} in the following senses: $a$) no ruler shorter then $N=6$ can exist that yields $K\frac{K-1}{2}=6$ distinct measures; and $b$) no further marks can be added to the ruler, without adding redundancy.
In general, these two distinct optimality criteria are defined as
\begin{itemize}
\item[$a$)] \emph{Length optimality}: Given a certain order $K$, the ruler's length $N$ is \emph{minimal};
\item[$b$)] \emph{Density optimality}: Given a certain length $N$, the ruler's order $K$ is \emph{maximal}.
\end{itemize}

The design of optimum Golomb rulers of higher orders is an NP-hard problem \cite{Apostolos2002, CottaCONSTRAINTS2007, Meyer2009}.
To illustrate the computational challenge involved, the Distributed.net project \cite{distributed2013}, which has the largest computing capacity in the world, has since the year 2000 dedicated a large share of its computing power to finding optimum Golomb rulers of various sizes.
The project took 4 years to compute the optimal Golomb ruler of order 24, and is expected to take 7 years to complete the search for the optimal oder-27 ruler!

%
%

\subsection{Genetic Algorithm to Design Orthogonal Golomb Rulers}
\label{Subsec:design_golomb_rulers}

Although a few systematic algorithms to generate Golomb rulers do exist \cite{Konstantinos2009},
none of the methods discovered so far are capable of outputting rulers adhering to a specific optimality criterion.
This, allied with the NP-hardness of the problem, makes efficient heuristic techniques the primary method to design Golomb rulers.
Indeed, the optimum rulers of orders 24 to 26 found by the Distributed.net were all obtained using heuristic methods \cite{Konstantinos2009, wiki:Golomb2013}.

Notice moreover that the optimality criteria described above are not necessarily sufficient to satisfy the needs of specific applications.
Due to the aforementioned reasons, heuristic techniques such as constraint programming \cite{Boese1994}, local search \cite{CottaCONSTRAINTS2007}, and evolutionary or genetic algorithms \cite{CarlosCotta2004, Ayari2010} are the standard approach to design Golomb rulers with specific features.

In the context of this article, our interest is to design \emph{orthogonal} Golomb rulers (so as to enable multipoint ranging), that also come as close as possible to satisfying the length and density of the optimality criteria described in subsection \ref{Subsec:Gol_Bas} (so as to optimise resources).
The orthogonality requirement adds the demand that rulers be designed out of a \underline{predefined} set of available integers $\mathcal{W}$, which to the best of our knowledge is an \emph{unsolved} problem.

In the next subsection we therefore describe a new genetic algorithm to design the required rulers.
The algorithm is a modified version of the technique first proposed in \cite{Soliday1995}, and inspired on the behaviour of wild animals that live in small groups, such as \emph{prides of lions}, and incorporate the following components.

\subsubsection{Representation}
\label{Subsubsec:Repre}

Following the framework proposed in \cite{Soliday1995}, Golomb rulers will be represented not by their marks, but by the differences of \emph{consecutive marks}.
That is, let $\mathcal{N} = \{n_1, n_2, \cdots,n_K\}$.
Then this set will be represented by $\mathcal{S} = \{s_1, \cdots, s_{K-1}\}$, where
\begin{equation}
\label{Eq:Repre}
s_i = n_{i+1} - n_i \quad \forall \quad i = \{1,\cdots,K-1\}.
\end{equation}

\subsubsection{Initial Population}
\label{Subsubsec:Ini_Pop}
In order to initialize the genetic algorithm an \emph{initial population} of segment sets is needed.
Let $s_{\max}$ be a design parameter describing the largest possible segment in the desired rulers, and consider the \emph{primary set} of segments $\mathcal{S}^* \triangleq  \{1,2,\cdots,s_{\max}\}$.
Then, each member $\mathcal{S}_p$ of the initial population is given by an $(K\!-\!1)$-truncation of a uniform random permutation of $\mathcal{S}^*$.

Notice that $s_{\max}$ must be larger then the order $K$ of the desired rulers, and that the larger the difference $s_{\max} - K$, the larger the degrees of freedom available to construct suitable rulers.

An initial population $\mathbb{P}$ of cardinality $p$ can then be defined as a set of $P$ non-equal segment sets $\mathcal{S}_p$, that is, $\mathbb{P} \triangleq \{\mathcal{S}_p\}_{p=1}^P$, with $\mathcal{S}_p \neq \mathcal{S}_q$ for all pairs $(p,q)$.

\subsubsection{Fitness Function}
\label{Subsubsec:Eva}

Once an initial population $\mathbb{P}$ is selected, each of the candidate rulers $\mathcal{S}_p$ are evaluated according to a fitness function designed to capture how closely the candidate ruler $\mathcal{S}_p$ approaches the prescribed features of the desired rulers.

Specifically, in the application of interest Golomb rulers must have length as small as possible for a given order (see optimality criteria in subsection \ref{Subsec:Gol_Bas}); must have all marks belonging to a certain set of admissible marks $\mathcal{W}$; and must have no repeated measures (by definition).

In order to define a suitable fitness function with basis on these criteria, let us denote the set of marks and the measure set corresponding to $\mathcal{S}_p$ respectively by $\mathcal{N}_p$ and $\mathcal{V}_p$.
Next, let $N_p$ and $F_p$ respectively denote the length and the minimum number of marks\footnote{Notice that in order to count $F_p$, all shifts of $\mathcal{N}_p$ within the range $[\min(\mathcal{W}),\max(\mathcal{W})]$ must be considered.} in $\mathcal{N}_p$ that are not in $\mathcal{W}$.
Finally, let $R_p$ be the number of repeated elements in $\mathcal{V}_p$.
Then, the fitness function is defined as
\begin{equation}
\label{Eq:Fitness}
f(\mathcal{S}_p) = N_p \times (R_p + F_p + 1).
\end{equation}

Notice that since randomly selected candidate rulers $\mathcal{S}_p$ are by construction suboptimal, $N_p \geq N$ for all $p$.
Furthermore, the sum $R_p+F_p$ is a non-negative integer, assuming the value $0$ only when no repetitions occur in
$\mathcal{V}_p$ and no marks outside $\mathcal{W}$ can be found in $\mathcal{N}_p$, \emph{simultaneously}.
In other words, the \emph{minimum} value of the fitness function is exactly $N$ and is achieved if and only if the respective candidate is indeed a Golomb ruler satisfying all the conditions required.

\subsubsection{Mutations}
\label{Subsubsec:Mut}
Although the fitness function has the desired property of being minimized only at optimum choices of $\mathcal{S}_p$, the underlying optimization procedure is not analytical, but combinatorial, due to the discreteness of the optimisation space (specifically, the space of all sets of segment sequences with $K-1$ elements).
Therefore, in order to optimize $f(\mathcal{S}_p)$ one needs to search the vicinity of $\mathcal{S}_p$, which is achieved by performing \emph{mutations} over the latter.

There are two distinct types of elementary mutations that can be considered: \emph{transmutation} and \emph{permutation}.
The first refers to the case where one element of $\mathcal{S}_p$ is changed to another value\footnote{Since a segment of length 1 is always required in a Golomb ruler \cite{Boese1994}, $s_i=1$ is never subjected to transmutation \cite{Boese1994}.}, while the second refers to a permutation between two segments.

Both types of mutation have similar effects on the all quantities $N_p$, $R_p$, and $F_p$.
But since a candidate sequence $\mathcal{S}_p$ is by definition already a Golomb ruler if $R_p = 0$, mutation is applied to $\mathcal{S}_p$ only if $R_p>0$.
And in that case, only one of the two types of elementary mutations is applied, randomly and with equal probability.

The elementary mutation operator will be hereafter denoted $\mathscr{M}(\cdot)$, and a version of $\mathcal{S}_p$ subjected to a single elementary mutation is denoted $\mathcal{S}^\dagger_p$ such we may write $\mathcal{S}^\dagger_p = \mathscr{M}(\mathcal{S}_p)$.

A sequence $\mathcal{S}_p$ is replaced by $\mathcal{S}^\dagger_p$ if and only if $f(\mathcal{S}^\dagger_p) < f(\mathcal{S}_p)$.
The mutation step is repeated for every $p$ until an improved replacement of $\mathcal{S}_p$ is found.
The mutation procedure is further iterated over the population $\mathbb{P}$ repeatedly until at least on candidate sequence $\mathcal{S}_p$ is Golomb, $R_p = 0$.
If no ruler can be found out of the initial population after a certain number of mutation iterations, the algorithm is restarted with an increased primary set  $\mathcal{S}^* \triangleq  \{1,2,\cdots,s_{\max}+1\}$.
This process is repeated until a mutated population $\mathbb{P}^\dagger$ is found, which contains at least one Golomb ruler.

\subsubsection{Selection}
\label{Subsubsec:Selection}
As a result of the mutation process described above, $\mathbb{P}^\dagger$ certainly contains one or more Golomb rulers.
Such rulers, however, may still violate the prescribed set of admissible marks $\mathcal{W}$ -- that is, may still have $F_p > 0$ -- and may not have the shortest length desired -- $i.e.$, $N_p > N$.

The optimized Golomb ruler will be obtained via the \emph{evolutionary} process to be described in the sequel, which in turn requires the classification the rulers in the population according to their function.
Specifically, the sequence $\mathcal{S}_p$ with $R_p = 0$ and the smallest score $f(\mathcal{S}_p)$ will be hereafter referred to as the \emph{dominant male} sequence and denoted $\mathcal{S}_{\male}$.
In other words, define $\mathbb{P}_{\!\male}^\dagger = \{\mathcal{S}_p| R_p = 0\}$, then
\begin{equation}
\mathcal{S}_{\male} = \{\mathcal{S}_p \in \mathbb{P}_{\male}^\dagger | f(\mathcal{S}_p) < f(\mathcal{S}_q)\;\forall\; q \neq p\}.
\end{equation}

In turn, all the other remaining sequences will be designated as \emph{female} sequences.
We shall therefore denote\footnote{Notice that this implies that ``male'' sequences in $\mathbb{P}_{\!\male}^\dagger$, but do not have the smallest score are thereafter relabelled ``female''.} $\mathbb{P}_{\female}^\dagger \triangleq \mathbb{P}^\dagger \setminus \mathcal{S}_{\male}$.

\subsubsection{Evolution}
\label{Subsubsec:Cross}
The evolution of the sequences occurs based on the Darwinian principle of variation via reproduction and selection by survival of the fittest.
Here, reproduction refers to the construction of new sequences via random crossover between the male sequence and any of the female ones, where crossover amounts to the swap of a block of adjacent ``genes'' from $\mathcal{S}_{\male}$ and $\mathcal{S}_{\female}$.

Let us denote the crossover operator as $\mathscr{C}(\cdot,\cdot)$, such that a child of $\mathcal{S}_{\male}$ and the $i$-th female $\mathcal{S}_{\female:i}$ in the population, generated via a single elementary crossover, can be described as $\mathscr{C}(\mathcal{S}_{\male},\mathcal{S}_{\female:i})$.
Then, the population evolves according to the following behaviour:
\begin{itemize}
\item The dominant male reproduces with all females generating the children $\mathscr{C}(\mathcal{S}_{\male},\mathcal{S}_{\female:i})$;
\item If there are any children with no repetition ($R_i=0$) and with fitness function lower then that of
$\mathcal{S}_{\male}$, then the child with the lowest score amongst those takes the place of the dominant male, that is
{\small
\begin{equation}
\mathcal{S}_{\male} \leftarrow \{\mathscr{C}(\mathcal{S}_{\male},\mathcal{S}_{\female:i}) | R_i = 0,\; f(\mathscr{C}(\mathcal{S}_{\male},\mathcal{S}_{\female:i})) < f(\mathcal{S}_{\male}) \;\text{and}\;f(\mathscr{C}(\mathcal{S}_{\female:i})) < f(\mathscr{C}(\mathcal{S}_{\female:j})) \;\forall\; j\neq  i \};
\end{equation}}
\item All other sequences are considered female, and out of original females and their children, only the best $P-1$ sequences, $i.e.$ the ones with the lowest scores, remains in $\mathbb{P}^\dagger$.
\end{itemize}

A pseudo-code of the genetic algorithm described above is given in Appendix \ref{Sec:PseudoCode}.
Due to the ``pride of lions'' evolutionary approach employed in the proposed algorithm, convergence to desired rulers is significantly faster then that achieved with the 
 ``giant octopus\footnote{It is known that both the female and male Pacific giant octopuses parish shortly after the hatching of their eggs \cite{howstuffworks:Octopus2012}.}'' approach taken in \cite{Soliday1995}, where both parent sequences are destroyed during the crossover process.

To illustrate the latter, consider the results shown in Table \ref{table:RelError}, which compares the 
average relative errors $\eta $ associated with Golomb rulers obtained with Soliday's algorithm \cite{Soliday1995} and the method proposed above, with
\begin{equation}
\eta \triangleq \mathbb{E}\left[\frac{N - N_\textup{opt}}{N_\textup{opt}}\right],
\end{equation}
where $N_\textup{opt}$ is the length of the shortest-possible (optimal) ruler with the same cardinality.

\begin{table}[H]
\centering
\caption{Comparison of Average Relative Error of Golomb Rulers}
\label{table:RelError}
\vspace{-3ex}
{\footnotesize
\hfill{}
\begin{tabular}{| c | c | c | c | c| }
\hline
$K$ & $N_\textup{opt}$ & {Soliday} \cite{Soliday1995} & {Proposed} ($P=2$) & {Proposed} ($P=4$)\\ \hline \hline
5 & 11 & 0.0\% & 0.0\% & 0.0\% \\ \hline
6 & 17 & 0.0\% & 0.0\% & 0.0\% \\ \hline
7 & 25 & 0.0\% & 0.0\% & 0.0\% \\ \hline
8 & 34 & 2.94\% & 0.0\% & 0.0\% \\ \hline
9 & 44 & 0.0\% & 4.6\% & 0.0\% \\ \hline
10 & 55 & 12.7\% & 12.7\% & 9.1\% \\ \hline
11 & 72 & 9.7\% & 11.1\% & 8.33\% \\ \hline
12 & 85 & 21.2\% & 16.5\% & 14.1\% \\ \hline
13 & 106 & 17.0\% & 17.0\% & 15.1\% \\ \hline
14 & 127 & 32.3\% & 23.6\% & 17.3\% \\ \hline
15 & 151 & 36.4\% & 26.5\% & 19.9\% \\ \hline
\end{tabular}}
\hfill{}
\vspace{-3ex}
\end{table}

It is found that even if the population considered during the evolution process is maintained to the minimum, replacing parents only by better offsprings tends to improve results as $K$ grows.
More importantly, a substantial and consistent improvement is achieved if $P>2$, such that the best (male) ruler can ``reproduce'' with multiple females.

Thanks to the modified fitness function (see equation \eqref{Eq:Fitness} compared to \cite[Eq. (4)]{Soliday1995}), which not only includes a direct term ($i.e.$, $F_p$) to account for the utilisation of forbidden marks, but also is only minimized when sequences are in fact Golomb rulers, the algorithm here proposed is capable of generating any desired number of orthogonal Golumb rulers, provided that $s_{\max}$ is sufficiently large.
This is achieved by subsequent executions of the algorithm, each time with $\mathcal{W}$ reduced by the marks of the rulers already generated.

There are, furthermore, two distinguished ways the resulting Golomb rulers can be grouped together.
One possibility is to group the rulers such that all have the same length $N$, even if with different  different number of marks.
This approach is motivated by the fact that the corresponding array-like vectors (see equation \eqref{Eq:Steerin_vect}) will have the same aperture, which in turn is directly related to the
accuracy of the corresponding distance estimation via superresolution algorithms.
This choice is referred to as Equivalent\footnote{As shall be demonstrated in Section \ref{Sec:CrlbPerf_Anal}, unequal Golomb rulers with the same $K$ and $N$, may still have different \acp{CRLB}.} Ranging Quality (ERQ) grouping.
Another possibility, however, is to group the Golomb rulers with the same cardinality $K$.
This grouping approach is motivated by the fact that, in the context hereby, each marker in the ruler corresponds to a measurement that is taken, and therefore is referred to as Fair Resource Allocation (FRA).

Examples of Golomb rulers obtained with the algorithm described above and grouped according to the ERQ and FRA criteria are listed in Table \ref{table:GolRul}.
It can be observed that, as desired, no two identical numbers can be found in two different rulers within the same group.
It follows that all the rulers of each group can be superimposed without interference and within a maximally compact span\footnote{If a conventional design were employed, the alternative would be to shift each ruler by length of the later!}.

To clarify, thanks to the rulers displayed in Table \ref{table:GolRul}, within a block of no more than $100$ cycles/frequencies, multipoint ranging between a source and 5 different anchors can be carried out by taking only $50$ \ac{ToA}/\ac{PDoA} measurements.
Furthermore, this can be achieved either with equivalent ranging quality using the group of ERQ rulers, or with fairly allocated resources using the group of FRA rulers, respectively.
\begin{table}[H]
\center
\caption{Examples of Golomb Rulers with ERQ and FRA Designs.}
\label{table:GolRul}
{\small
\hfill{}
\begin{tabular}{|c|l|c:c||c|l|c:c|}
\hline
$K$ &\multicolumn{1}{c|}{{Equal Ranging Quality}}& $N$& $M$ & $K$ &\multicolumn{1}{c|}{{Fair Resource Allocation}}& $N$& $M$ \\
\hline
\hline 
9 & 0,1,7,10,30,41,45,63,87 			&87 &36 & 10 & 0,1,16,21,24,49,63,75,81,85   &85 &45 \\
\hline 
9 & 2,3,6,32,37,49,56,76,89 			&87 &36 & 10 & 2,3,11,32,45,56,60,72,78,92   &90 &45 \\
\hline 
10 & 4,5,16,20,33,42,52,66,73,91 		&87 &45 & 10 & 5,9,15,29,42,51,68,80,91,96   &91 &45 \\
\hline 
11 & 8,9,18,21,38,46,53,72,77,93,95 	&87 &55 & 10 & 6,13,17,19,33,43,61,62,84,93  &87 &45 \\
\hline 
11 & 12,13,17,25,31,47,68,70,79,96,99 	&87 &55 & 10 & 12,14,22,27,28,46,66,73,77,94 &82 &45 \\
\hline 
\end{tabular}}
\hfill{}
\vspace{-3ex}
\end{table}

\section{Error Analysis and Comparisons}
\label{Sec:CrlbPerf_Anal}

In this section we analyse the performance of the multipoint ranging approach described above, both with \ac{PDoA} and \ac{ToA} measurements.
To this end, we first derive the Fisher Information Matrices and associated Cramer-Rao Lower Bounds (CRLB) corresponding to the algorithms and later offer comparisons with simulated results.
Since related material on \ac{ToA} can be found more easily \cite{Tao2008,Kaune2012}, we shall consider first the \ac{PDoA} case and offer only a synthesis of the \ac{ToA} counterpart.
%
%

\vspace{-2ex}
\subsection{Phase-Difference of Arrival}
\label{Subsubsection:PDoA}

Start by recognising that phase difference measurements subject to errors are circular random variables.
The Central Limit Theorem (CLT) over circular domains establishes that the most entropic ($i.e.$, least assuming) model for circular variables with known mean and variance is the von Mises or Tikhonov distribution \cite{Abreu2008}.
We assume, therefore, that phase measurements are modeled as
\begin{equation}
\label{Eq:Tikhonov_angles}
\hat{\Delta}\varphi  \sim P_\mathcal{T}(x;\Delta\varphi,\kappa)
\end{equation}
with
\begin{equation}
\label{Eq:Tikhonov_pdf1}
P_\mathcal{T}(x;\Delta\varphi,\kappa) \triangleq \dfrac{1}{2\pi I_0(\kappa)}\cdot \exp(\kappa\cos(x - \Delta \varphi)),\quad -\pi \leq x \leq \pi,
\end{equation}
where $I_n(\kappa)$ is the $n$-th order modified Bessel function of the first kind and $\kappa$ is a shape parameter which in the case of phase estimates is in fact given by the signal-to-noise-ratio (SNR) of input signals, and that relates to the error variance by
\begin{equation}
\sigma^2_{\Delta\varphi} = {1 - \frac{I_1(\kappa)}{I_0(\kappa)}}
\mathrel{\mathop{\kern0pt\xrightarrow{\hspace*{3em}}}\limits_{\kappa >> 1}} {\frac{2}{2\kappa + 1}} \approx \frac{2}{\kappa}.
\end{equation}

Consider then that a set of $K$ independent measurements $\{\Delta \varphi_k\}_{k\in\mathcal{N}}$ is collected according to a Golomb ruler $\mathcal{N}$, such that the samples can be expanded into and augmented set of $M$ samples $\{\Delta \varphi_m\}_{m\in\mathcal{V}}$, with
\begin{equation}
\label{Eq:DeltaPhiGolomb}
\Delta \varphi_m = \Delta \varphi_k - \Delta \varphi_\ell = \omega_d (k - \ell) = \omega_d \nu_m, \quad \mbox{for}\quad k > \ell\quad \mbox{and}\quad (k,\ell)\to m,
\end{equation}
where each index $m$ corresponds to a pair $(k,\ell)$ with $k>\ell$ with ascending differences\footnote{Notice that this is ensured without ambiguity thanks to the fact that $\mathcal{N}$ is a Golomb ruler.}, and we commit a slight abuse of notation compared to equation \eqref{Eq:PhaseDifferenceDistance}, since $\nu_m$ is a positive integer obtained from a the difference $k-\ell$, such that $\nu_m \neq m$.

At this point it is worthy of mention that although the expanded samples $\Delta \varphi_m$ are actually differences of phase differences, these quantities not only preserve the linear relationship with the parameter of interest but also their independence.
As a result of the double-differences, however, the SNR of $\Delta \varphi_m$ from equation \eqref{Eq:DeltaPhiGolomb} is twice that of $\Delta \varphi_k$ from equation \eqref{Eq:PhaseDifferenceDistance}.
In light of the asymptotic relationship being twice as large, it follows that the shape parameter $\kappa$ associated with $\Delta \varphi_m$'s are twice as small.

Using the model above, and incorporating the optimised sampling via Golomb ruler, the likelihood function associated with $M$ independent measurements as per equation \eqref{Eq:PhaseDifferenceDistance} becomes,
\begin{equation}
\label{Eq:Likelihood}
L_\mathcal{T}(\hat{d};\Delta f,\kappa) = \prod_{m = 1}^{M} P_\mathcal{T}(x;\Delta\varphi_m ,\kappa) = \displaystyle \dfrac{1}{(2\pi I_0(\kappa/2))^M} \prod_{m = 1}^{M}
\exp\left[\frac{\kappa}{2}\cos\left(\frac{4\pi\Delta f}{c} \nu_m\cdot(\hat{d} - d)\right)\right],
\end{equation}
where $\nu_m \in \mathcal{V}$ and we have slightly modified the notation in order to emphasise the quantity and parameter of interest $d$.

For future convenience, let us define $\alpha = \frac{4\pi\Delta f}{c}$. Then the associated log-likelihood function is
\begin{equation}
\label{Eq:Loglikelihood}
\ln L_\mathcal{T}(\hat{d};\Delta f,\kappa) = -M \ln {2\pi I_0(\kappa/2)} + \frac{\kappa}{2} \sum_{m=1}^{M}
\cos\left(\alpha\cdot \nu_m\cdot (\hat{d} - d)\right),
\end{equation}
and its Hessian becomes
\begin{equation}
\label{Eq:Hess_loglikehood}
\dfrac{\partial ^2 \ln L_\mathcal{T}(\hat{d};\Delta f,\kappa)}{\partial \hat{d}^2} =
- \frac{\alpha^2\kappa}{2}\sum_{m=1}^{M} \nu_m^2\cos\left(\alpha\cdot \nu_m\cdot(\hat{d} - d)\right).
\end{equation}

The Fisher Information is the negated expectation of the Hessian, thus,
\begin{equation}
J(\mathcal{V};\Delta f,\kappa) = -\mathbb{E} \left[\dfrac{\partial ^2 \ln L_\mathcal{T}(\hat{d};\Delta f,\kappa)}{\partial \hat{d}^2}\right] = \frac{\alpha^2\kappa}{2}\sum_{m=1}^{M} \nu_m^2\mathbb{E}\left[\cos\left(\alpha\cdot \nu_m\cdot(\hat{d} - d)\right)\right],
\end{equation}
where the notation alludes to the fact that the key input determining the Fisher Information is the set of measures $\mathcal{V}=\{\nu_1,\cdots,\nu_M\}$.

Next, recognise that each term $\alpha\cdot \nu_m\cdot(\hat{d} - d)$ is in fact a centralized circular variate with the same distribution $P_\mathcal{T}(x;0,\kappa/2)$, regardless of $m$.
Then, substituting $\alpha\cdot \nu_m\cdot(\hat{d} - d)$ with $\theta$, we obtain
\begin{align}
J(\mathcal{V};\Delta f,\kappa) & = \frac{\alpha^2\kappa}{2}
	\sum_{m=1}^{M} \nu_m^2\mathbb{E}\left[\cos\theta\right] = \frac{\alpha^2\kappa}{2} \sum_{m=1}^{M} \dfrac{\nu_m^2}{I_0(\kappa/2)}\underbrace{\frac{1}{\pi}\int\limits_{0}^{\pi} \cos\theta \exp\left(\frac{\kappa}{2}\cos \theta\right) \,\d \theta}_{I_1(\kappa/2)}\nonumber\\[-4ex]
  & = \frac{\alpha^2\kappa}{2} \dfrac{I_1(\kappa/2)}{I_0(\kappa/2)}\sum_{m=1}^{M} \nu_m^2,
\end{align}
where the integration limits in the integral above follow from evenness of the function $\cos(\theta)\exp(\frac{\kappa}{2}\cos\theta)$, and the last equality results from the integral solution found in \cite[Eq. 9.6.19, pp. 376]{MyListOfPapers:Abramowitz1965}.

Since the above Fisher Information is a scalar, the CRLB is obtained directly by taking its inverse, $i.e.$,
\begin{equation}
\label{Eq:Crlbrang1}
\text{CRLB}_\textup{PDoA}(\mathcal{V};\Delta f,\kappa) = \frac{1}{J(\mathcal{V};\Delta f,\kappa)}.
\end{equation}

Before proceeding to the \ac{ToA} case, some discussion on the analytical results offered above are in order.
First, let us emphasise that given a set of phase difference measurements $\{\Delta\varphi_{n_k}\}_{k=1}^{K}$, with $n_k\in \mathcal{N}$, one always has the \emph{option} of either exploit the properties of the Golomb ruler $\mathcal{N}$ and expand to a set of measurements $\{\Delta\varphi_{\nu_m}\}_{m=1}^{M}$, or not.
In case such option is \emph{not} adopted, the associated Fisher Information and CRLB can obviously be obtained exactly as done above, but with $\kappa$ replacing $\kappa/2$ and $\mathcal{N}$ replacing $\mathcal{V}$.
That is,
\begin{equation}
\label{Eq:CrlbrangNotExpanded}
J(\mathcal{N};\Delta f,\kappa) = {\alpha^2\kappa} \dfrac{I_1(\kappa)}{I_0(\kappa)}\sum_{k=1}^{K} n_k^2 \quad \Longleftrightarrow \quad \text{CRLB}_\textup{PDoA}(\mathcal{N};\Delta f,\kappa) = \frac{1}{J(\mathcal{N};\Delta f,\kappa)}.
\end{equation}

Comparing these expressions, it can be readily seen that the choice of adopting the Golomb approach on the one hand subjects the resulting double-phase-differences to twice the noise, but on the other hand expands the number terms in the summation.
In principle, the optimum choice between these options therefore depends on the ruler $\mathcal{N}$ and its order $K$, and the associated $\mathcal{V}$ and $M$, as well as $\kappa$.
As can be shown in Figure \ref{Fig:crlb_kappa}, for instance, the ruler $\mathcal{N}=\{0,1,4,6\}$ yields superior results compared to its associated measure set $\mathcal{V}=\{1,2,3,4,5,6\}$, because the loss of 3dB (implied by $\kappa \to \kappa/2$) incurred by the latter is not compensated by the increase gained in the sum of squares achieved by using $\mathcal{V}$ instead of $\mathcal{N}$.

For larger rulers, however, the advantage of expanding the rulers quickly becomes significant, thanks to the geometric increase of $M$ with respect to $K$.
A ruler of order 6, $e.g.$, $\mathcal{N}=\{0,1,4,10,12,17\}$, already achieves better performance expanded into
$\mathcal{V}=\{1,\cdots,17\}$ then otherwise, for $\sigma_{\Delta_\varphi} \leq 0.22$.
Likewise, the expanded version of the order-10 ruler $\mathcal{N}=\{0,1,16,21,24,49,63,75,81,85\}$ is superior up to $\sigma_{\Delta_\varphi} \leq 0.65$ -- which incidentally defines essentially the entire range of interest --
and finally the expanded ruler of order-20 is always superior, for any $\sigma_{\Delta_\varphi}$.
In summary, it can be said that applying the Golomb expansion leads to superior results, as long as the ruler is large enough and $\sigma_{\Delta_\varphi}$ is in the region of interest.
%
%

\subsection{Time of Arrival}
\label{Subsubsection:ToA}

Due to the similarity of the \ac{ToA} and \ac{PDoA} ranging models described in Section \ref{Sec:Prelim}, the Fisher Information and CRLB for \ac{ToA}-ranging with Golomb rulers are very similar to those given above for the \ac{PDoA} case.
For the sake of brevity, we therefore offer here only a succinct derivation.

Assuming that the error on the time of arrival estimates are Gaussian-distributed, were have
\begin{equation}
\label{Eq:time_arrival}
\hat{\Delta}\tau \sim P_\mathcal{G}(x;\Delta\tau,\sigma_{\Delta\tau}^2) =  \dfrac{1}{\sqrt{2\pi} \sigma_{\Delta\tau}}
\exp\left(-\frac{(x - \Delta\tau)^2}{2\sigma_{\Delta\tau}^2}\right),
\end{equation}
such that the likelihood function, the log-likelihood function, its Hessian and the Fisher Information, considering already the expansion $\mathcal{N}\to\mathcal{V}\Rightarrow \sigma_{\Delta\tau}^2 \to 2\sigma_{\Delta\tau}^2$ and emphasising the quantities of interest, becomes
\begin{eqnarray}
\label{Eq:LikelihoodGaussian}
&L_\mathcal{G}(\hat{d};\sigma^2_{\Delta\tau}) = \prod_{m = 1}^{M} P_\mathcal{G}(\hat{d};\Delta\tau_m,2\sigma_{\Delta\tau}^2) =
\displaystyle \dfrac{1}{(4\pi\sigma_{\Delta\tau}^2)^{M/2}} \prod_{m = 1}^{M}
\exp\left(-\frac{\nu_m^2}{c^2\sigma_{\Delta\tau}^2}(\hat{d} - d)^2\right),&\nonumber\\
&\ln L_\mathcal{G}(\hat{d};\sigma^2_{\Delta\tau}) = -\frac{M}{2} \ln 4\pi\sigma^2_{\Delta\tau} - \dfrac{1}{c^2\sigma^2_{\Delta\tau}}\displaystyle\sum_{m=1}^{M} \nu_m^2 (\hat{d}-d)^2,&\nonumber\\
&\dfrac{\partial ^2 \ln L_\mathcal{G}(\hat{d};\sigma^2_{\Delta\tau})}{\partial \hat{d}^2} = -\dfrac{2}{c^2\sigma^2_{\Delta\tau}}\displaystyle\sum_{m=1}^{M} \nu_m^2 \quad\Longrightarrow\quad J(\mathcal{V};\sigma^2_{\Delta \tau}) = \dfrac{2}{c^2\sigma^2_{\Delta\tau}}\displaystyle\sum_{m=1}^{M} \nu_m^2.& \label{Eq:FisherToAGolomb}
\end{eqnarray}

As discussed above, if the measurements taken according to the Golomb markers are, however, used without taking their differences, the associated noise process has half the variance such that
\begin{equation}
J(\mathcal{N};\sigma^2_{\Delta \tau}) = \dfrac{4}{c^2\sigma^2_{\Delta\tau}}\displaystyle\sum_{m=1}^{M} \nu_m^2.\label{Eq:FisherToANoGolomb}
\end{equation}

\subsection{Simulations and Comparison Results}
\label{Sec:Results}

Let us finally study the performance of the proposed multipoint ranging technique by means of simulations and comparisons with the corresponding CRLBs derived above.
For the sake of brevity, we will consider only \ac{PDoA} ranging as all results obtained with the \ac{ToA} approach are equivalent.

\begin{figure}[H]
\centering
\subfigure[{\footnotesize As function of $K$, for different $\sigma_{\Delta\varphi}$.}]%
{\includegraphics[height=0.4\textheight]{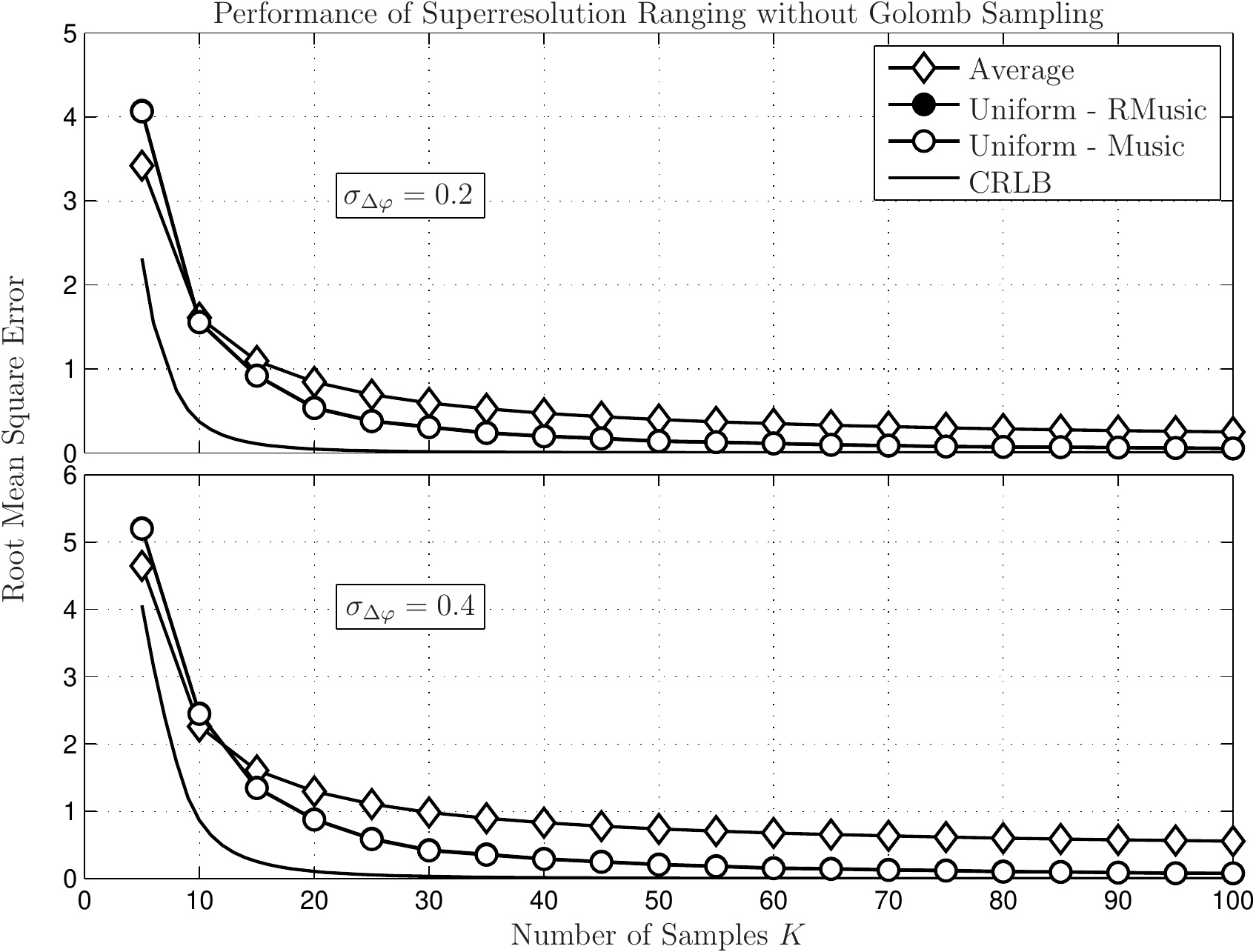}
\label{Fig:Stra_musrmus_freq}}
\vspace{1ex}
\subfigure[{\footnotesize As function of $\sigma_{\Delta\varphi}$, for different $K$.}]%
{\includegraphics[height=0.4\textheight]{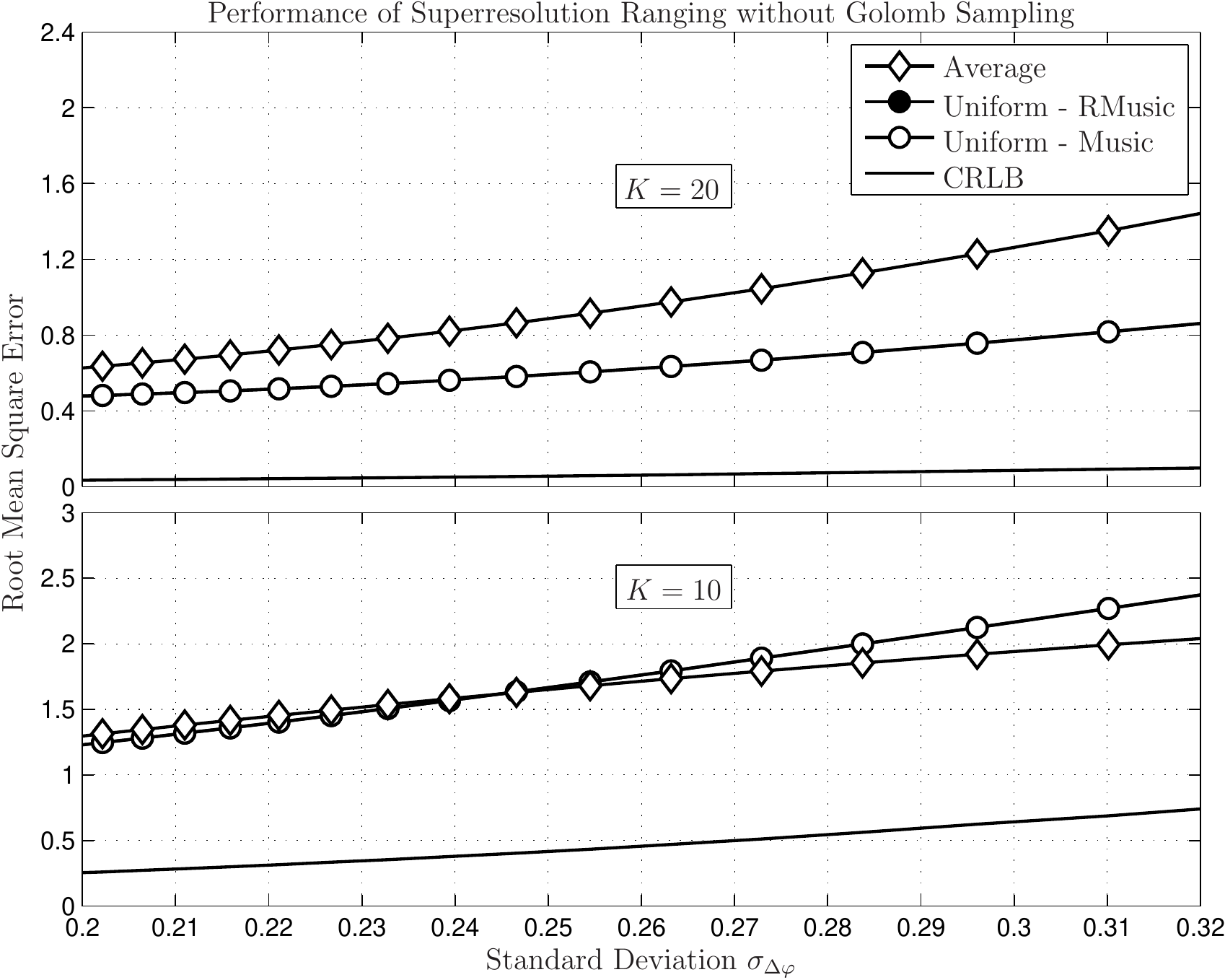}
\label{Fig:Stra_musrmus_kappa}}
\vspace{-2ex}
\caption{Performance of superresolution and average-based ranging algorithms as a function of the sample set sizes $K$ and the phase-estimate noise variances $\sigma_{\Delta\varphi}$, \emph{without} Golomb-optimized sampling.}
\label{Fig:Stra_musrmus_freqkappa} 
\end{figure}

First, consider Figure \ref{Fig:Stra_musrmus_freqkappa}, where the performances of two classic superresolution algorithms -- namely the Music and Root Music algorithms of briefly described in Subsection \ref{Subsec:S_PDB_R} -- are compared against the CRLB derived in Subsection \ref{Subsubsection:PDoA}.
Plots are shown both as a function of $K$ for various $\sigma_{\Delta\varphi}$ and vice-versa, and for the sake of having a practical reference, we include also results obtained by simply averaging the distance estimates corresponding to all independent samples.

We emphasise that in this figure \underline{no Golomb ruler is used}.
Instead, a sequence of $K$ consecutive samples is collected for each range estimate, as typically assumed in existing work \cite{Chien-Sheng2012, Junyang2012}.

One fact learned from these plots -- and is particularly visible in Figure \ref{Fig:Stra_musrmus_freq} -- is that without the efficient use of samples made possible by the Golomb ruler approach here proposed, superresolution algorithms require a large number of samples in order to reach the CRLB, which is a problem since energy consumption and latency are directly related to the number of samples collected. 

Another fact of relevance that can be learned, however, is that although supperresolution methods do improve on a ``naive'' average-based estimator, that gain in itself is not that significant unless the number of samples $K$ is rather large.
This is highlighted in Figure \ref{Fig:Stra_musrmus_kappa}, where it is seen that with $K=10$, the simple average-based algorithm has essentially the same performance of MUSIC.

\begin{figure}[H]
\centering
\subfigure[{\footnotesize As function of $K$, for different $\sigma_{\Delta\varphi}$.}]%
{\includegraphics[height=0.4\textheight]{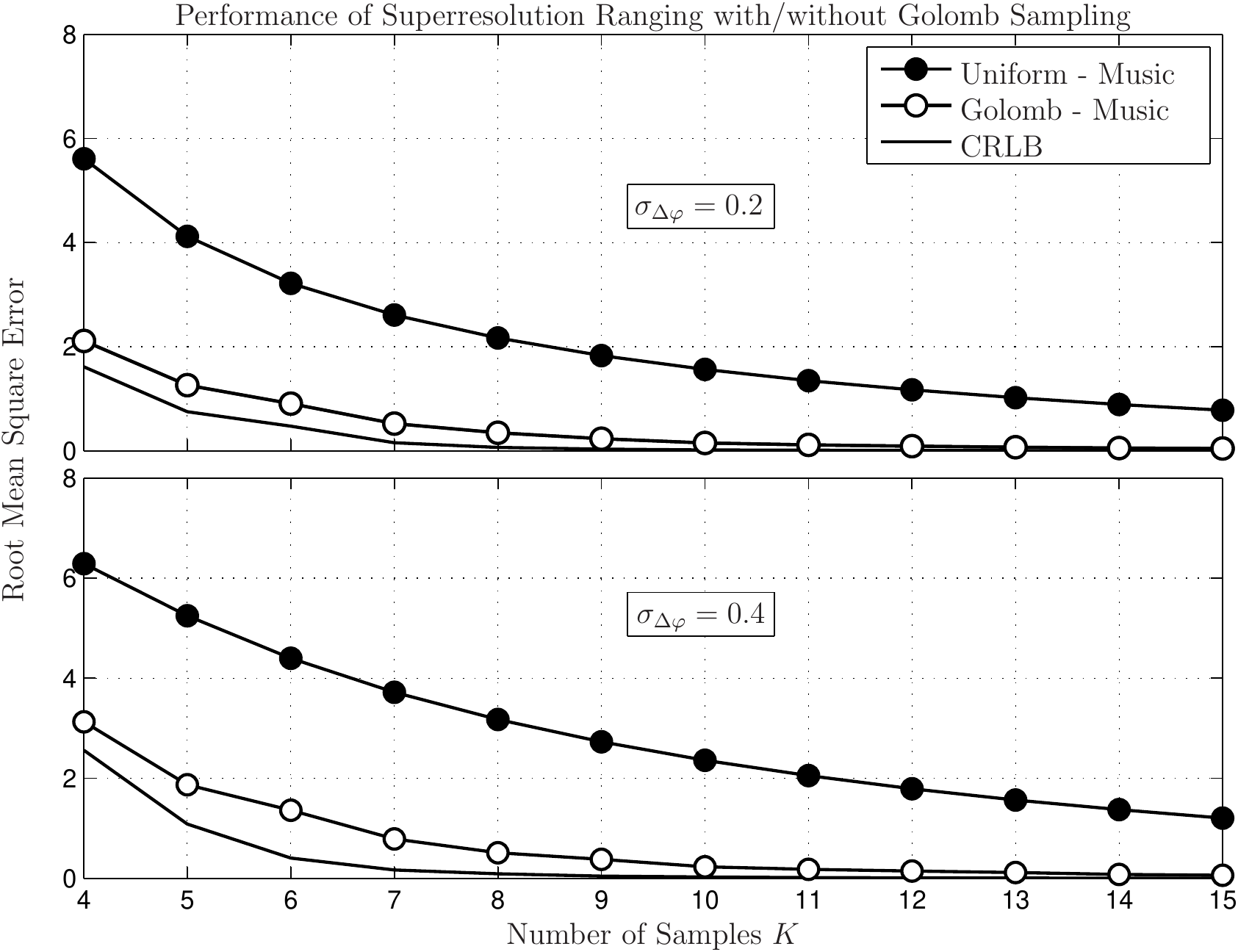}
\label{Fig:Gol_rmus_freq}}
\subfigure[{\footnotesize As function of $\sigma_{\Delta\varphi}$, for different $K$.}]%
{\includegraphics[height=0.4\textheight]{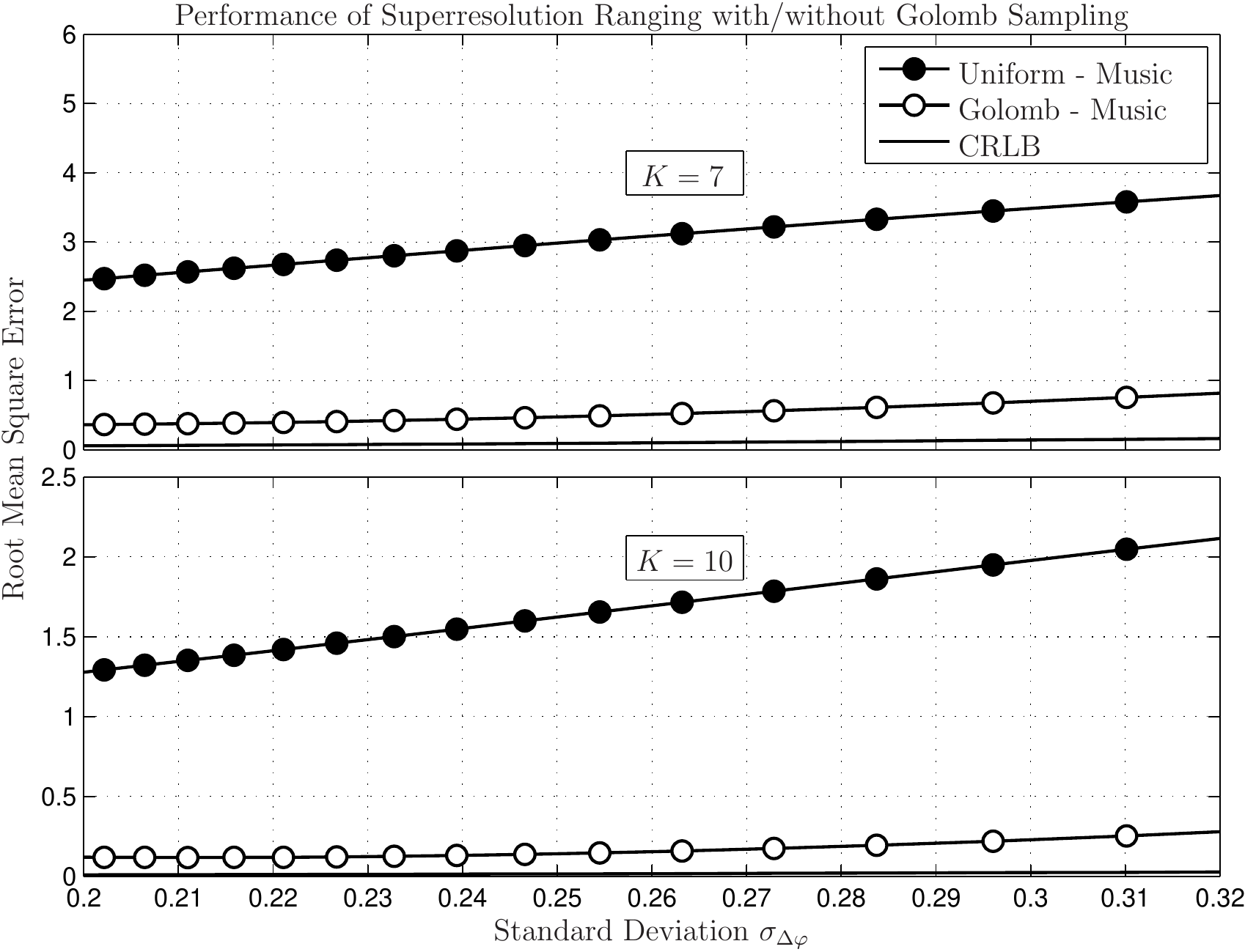}
\label{Fig:Gol_rmus_kappa} }
\vspace{-1ex}
\caption{Performance of superresolution ranging algorithms as a function of the sample set sizes $K$ and the phase-estimate noise variances $\sigma_{\Delta\varphi}$, both \emph{with} and \emph{without} Golomb-optimized sampling.}\label{Fig:Gol_rmus_freqkappa} 
\end{figure}

The results above emphasize the significance of our contribution, by demonstrating that the efficient utilisation of samples is fundamental to reap from superresolution algorithms their true potential performance.
This is further illustrated in Figure \ref{Fig:Gol_rmus_freq}, where it can be seen that thanks to the Golomb sampling superresolution algorithms with a relatively small number of samples come much closer to the CRLB.

Considered in coordination with the results of Figure \ref{Fig:crlb_kappa}, it can be generally said that a Golomb-optimized scheme with a total of 10 samples, taken at frequencies corresponding to an accordingly Golomb ruler $\mathcal{N}$ expanded into the associated measure set $\mathcal{V}$, followed by MUSIC estimation is an excellent choice for \ac{PDoA} ranging.

In fact, as illustrated by Table \ref{table:GolRul}, such a choice also allows for an easy design of various orthogonal Golomb rulers, such that multipoint ranging can be efficiently performed.
But since in this case a choice needs to be made between the ERQ and FRA ruler allocation approaches, a fair question to ask in this context is what are the performances of corresponding choices.

\begin{figure}[H]
\centering
\includegraphics[height=0.4\textheight]{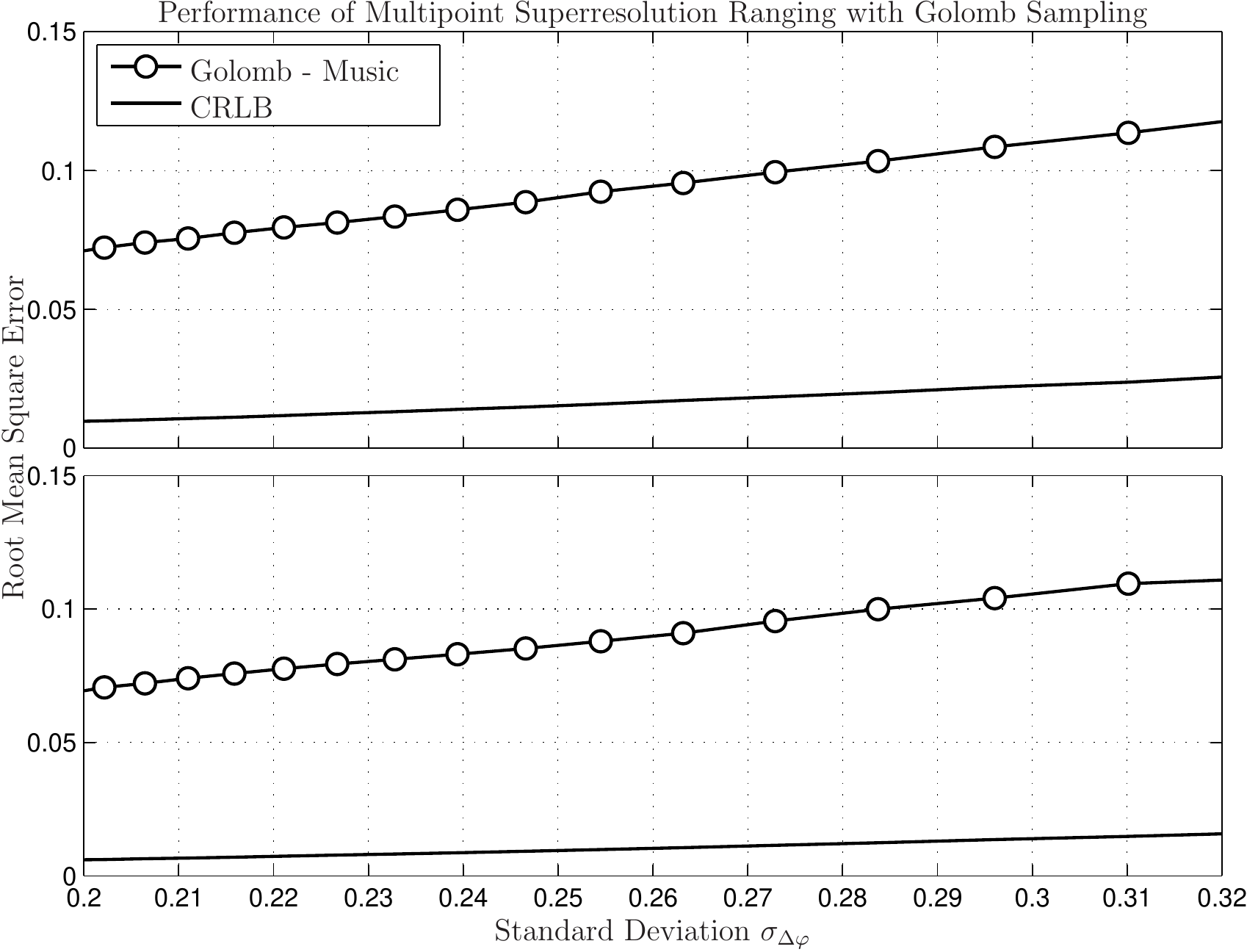}
\vspace{-3ex}
\caption{Performance of Golomb-optimized superresolution ranging with ERQ and FRA ruler allocation approaches.}
\label{Fig:Multipoint} 
\end{figure}

This is addressed in Figure \ref{Fig:Multipoint}, where the average performances of an ERQ and an FRA multipoint ranging schemes employing the rulers shown in  Table \ref{table:GolRul} are compared against corresponding CRLBs.
The figure shows that in fact both approaches have similar performances relative to one another and relative to the CRLBs.


\section{Conclusions}
\label{Sec:Conclusion}
We offered an efficient and accurate solution to the multipoint ranging problem, based on an adaptation of superresolution techniques, with optimised sampling.
Specifically, using as examples the specific cases of \ac{ToA} and \ac{PDoA}, unified under the same mathematical framework, we constructed a variation of the MUSIC and Root-MUSIC algorithm to perform distance estimation over sparse sample sets determined by Golomb rulers.
The design of the mutually orthogonal sets of Golomb rulers required by the proposed method -- a problem that founds no solution in current literature -- was shown to be achievable via a new genetic algorithm, which was also shown to outperform the best known alternative when used to generate optimal rulers.
A \acp{CRLB} analysis of the overall optimised multipoint ranging solution was performed, which compared to simulated results quantified the substantial gains achieved by the proposed technique.


\section{Acknowledgements}
This work has been performed within the framework FP7 European Union Project BUTLER (grant no. 287901).


\newpage
\appendices
\section{}
\label{Sec:PseudoCode}
\vspace{-1ex}
\begin{algorithm}[H]
{\small
\caption{- Golomb Ruler Generation Algorithm}\label{gen_algo}
\begin{algorithmic}
\State $\mathcal{W} \longleftarrow$ Set of forbidden marks (given)
\vspace{-1ex}
\State $K \longleftarrow$ Desired order of the ruler (given)
\vspace{-1ex}
\State $C \longleftarrow$ Maximum number of mutations (given)
\vspace{-1ex}
\State $G \longleftarrow$ Maximum number of generations (given)
\vspace{-1ex}
\State $s_{\textup{max}}:= K-1$
\vspace{-1ex}
\While{$\nexists\; \mathcal{S}_p | f(\mathcal{S}_p) = K\frac{(K-1)}{2}$}
\vspace{-1ex}
 \State $\mathcal{S}^* := \{1,2,\cdots,s_{\textup{max}}\}$
 \vspace{-1ex}
 \For{$p:= 1 \to P$}
  \vspace{-1ex}
  \State count $\leftarrow 0$ 
  \vspace{-1.5ex}
  \State $\mathcal{S}_p \leftarrow$ randomly select $K-1$ elements of $\mathcal{S}^*$
  \vspace{-1.5ex}
  \State $\mathcal{S}_p \leftarrow$ randomly permute the elements of $\mathcal{S}_p$
  %
  \vspace{-1ex}
  \While{count $<=$ $C$}
   \vspace{-1ex}
   \State $\mathcal{S}^\dagger_p = \mathscr{M}(\mathcal{S}_p)$.
   \vspace{-1ex}
   \If{$f(\mathcal{S}^\dagger_p) < f(\mathcal{S}_p)$}
    \vspace{-1ex}
    \State $S_p \leftarrow \mathcal{S}^\dagger_p$
    \vspace{1ex}
   \EndIf{\bf end if} 
   \vspace{-1.5ex}
   \State count $\leftarrow$ count + $1$
 \vspace{1ex}
 \EndWhile{\bf end while} 
 \vspace{1ex}
 \EndFor{\bf end for} 
 \vspace{-1.5ex}
 \State $\mathbb{P}^\dagger = \{\mathcal{S}_p\}^{P}_{p=1}$
 %
 \vspace{-1ex}
 \If{$\nexists\; \mathcal{S}_p | R_p=0$}
  \vspace{-1ex}
  \State $s_{\textup{max}} \leftarrow s_{\textup{max}}+1$
  \vspace{-1.5ex}
  \State restart
  \vspace{-1.5ex}
 \Else
 \vspace{-1.5ex}
 %
  \State $\mathbb{P}_{\!\male}^\dagger \leftarrow \{\mathcal{S}_p| R_p = 0\}$
  \vspace{-1ex}
  \State $\mathcal{S}_{\male} =
    \{\mathcal{S}_p \in \mathbb{P}_{\male}^\dagger | f(\mathcal{S}_p) <f(\mathcal{S}_q)\;\forall\; q \neq p\}$
  \vspace{-1ex}
  \State $\mathbb{P}_{\female}^\dagger \leftarrow \mathbb{P}^\dagger \setminus \mathcal{S}_{\male}$
  \vspace{1ex}
 \EndIf{\bf end if}
 \vspace{-1.5ex}
 %
 \State count $\leftarrow$ 0
 \vspace{-1.5ex}
 \While{count $<$ G {\bf or} $\exists\; \mathcal{S}_p | f(\mathcal{S}_p) = K\frac{(K-1)}{2}$}
  \vspace{-1ex}
  \For{$p := 1 \to P-1$}
   \vspace{-1ex}
   \State $\mathcal{S}^\dagger_p \leftarrow \mathscr{C}(\mathcal{S}_p,\mathcal{S}_{\male})$
   \vspace{-1ex}
   \If{$f(\mathcal{S}^\dagger_p) < f(\mathcal{S}_p)$}
    \vspace{-1ex}
    \State $\mathcal{S}_p \leftarrow \mathcal{S}^\dagger_p$
    \vspace{-1ex}
   \ElsIf{$f(\mathcal{S}^\dagger_p) < f(\mathcal{S}_{\male})$}
    \vspace{-1ex}
    \State $\mathcal{S}_{\male} \leftarrow \mathcal{S}^\dagger_p$
    \vspace{1ex}
   \EndIf{\bf end if}
   \vspace{1ex}
  \EndFor{\bf end for}
  \vspace{1ex}
 \EndWhile{\bf end while}
 \vspace{1ex}
\EndWhile{\bf end while}
\vspace{-1ex}
%
\Function{Fitness Function}{}
 \vspace{-1ex}
 \State $N_p \leftarrow$ length of $\mathcal{N}_p$ associated to $\mathcal{S}_p$ (input)
 \vspace{-1ex}
 \State $R_p \leftarrow$ number of repeated elements in $\mathcal{S}_p$ (input)
 \vspace{-1ex}
 \State $F_p \leftarrow$ number of forbidden marks in $\mathcal{N}_p$ (input)
 \vspace{-1ex} 
 \State $f(\mathcal{S}_p) \leftarrow N_p \times (R_p + F_p + 1)$.
 \vspace{-1ex}
 \State \Return $f(\mathcal{S}_p)$
\vspace{0.5ex}
\EndFunction{\bf end function}
\end{algorithmic}
}
\end{algorithm}

\newpage





\newpage

\newpage
\newpage
%

\newpage

\newpage

\newpage


\end{document}